\documentclass[11 pt, reqno]{article}

\usepackage{amsmath,amsfonts,amssymb,amsthm}
\usepackage[colorlinks=true,hyperindex=true]{hyperref}
\usepackage{cancel}
\usepackage{mathabx} 
\usepackage{comment}
\usepackage{enumitem}
\usepackage{cite}
\usepackage{showlabels}

\usepackage{chngcntr}
\counterwithin*{equation}{section}

\usepackage[a4paper, left=2cm, right=2cm, top=2 cm, bottom= 2 cm]{geometry}
\usepackage{color}
\usepackage{fancyhdr}
\usepackage{latexsym}

\newtheorem{ccounter}{ccounter}[section]
\newtheorem{thm}[ccounter]{Theorem}
\newtheorem{lem}[ccounter]{Lemma}
\newtheorem{cor}[ccounter]{Corollary}
\newtheorem{defn}[ccounter]{Definition}
\newtheorem{prop}[ccounter]{Proposition}
\newtheorem{ass}[ccounter]{Assumption}
\newtheorem{ex}[ccounter]{Example}
\newtheorem{definition}[ccounter]{Definition}

\def\bet{\begin{thm}}
\def\eet{\end{thm}}
\def\bel{\begin{lem}}
\def\eel{\end{lem}}
\def\bas{\begin{ass}}
\def\eas{\end{ass}}
\def\bec{\begin{cor}}
\def\eec{\end{cor}}
\def\bed{\begin{defn}}
\def\eed{\end{defn}}
\def\bep{\begin{prop}}
\def\eep{\end{prop}}
\def\beq{\begin{equation}}
\def\eeq{\end{equation}}
\def\proof{\noindent {\bf Proof.}\ \ }
\def\bea{\begin{equation*}}
\def\eea{\end{equation*}}
\def\tr{\mathrm{tr}}
\def\bex{\begin{ex}}
\def\eex{\end{ex}}
\def\remark{\noindent{\bf Remark. }}
\def\ulj{\underline{j}}
\def\uli{\underline{i}}
\def\P{\mathcal{P}}
\def\chit{\chi^{(2)}}
\def\rhosc{\rho_{\mathrm{sc}}}
\def\msc{m_{\mathrm{sc}}}

\def\rr{\mathbb{R}}

\def\cc{\mathbb{C}}
\def\1{\boldsymbol{1}}
\def\Im{\mathrm{Im}}
\def\Re{\mathrm{Re}}
\def\e{\mathrm{e}}
\def\i{\mathrm{i}}
\def\del{\partial}
\def\d{\mathrm{d}}
\def\eps{\varepsilon}
\renewcommand\leq\varleq
\renewcommand\geq\vargeq
\def\ee{\mathrm{E}}

\def\F{\mathcal{F}}
\def\O{\mathcal{O}}

\def\ee{\mathbb{E}}

\def\pp{\mathbb{P}}
\def\cb{c_\beta}
\def\tilm{\tilde{m}}
\def\N{\mathcal{N}}
\def\G{\mathcal{G}}
\def\F{\mathcal{F}}

\def\A{\mathcal{A}}

\begin{document}
\title{Deformed GOE}


\begin{table}
\centering

\begin{tabular}{c}
\Large{\bf Fluctuations of the overlap at low temperature}\\
 \Large{\bf in the 2-spin spherical SK model}
\\
\\
\end{tabular}
\begin{tabular}{c c c}
Benjamin Landon & & Philippe Sosoe \\
 & & \\
 \small{Massachusetts Institute of Technology} & & \small{Cornell University} \\
  & & \\
\end{tabular}
\\
\begin{tabular}{c}
\multicolumn{1}{c}{\today}\\
\\
\end{tabular}

\begin{tabular}{p{15 cm}}
\small{{\bf Abstract:} We describe the fluctuations of the overlap between two replicas in the 2-spin spherical SK model about  its limiting value in the low temperature phase. We show that the  fluctuations are of order $N^{-1/3}$ and are given by a simple, explicit function of the eigenvalues of a matrix from the Gaussian Orthogonal Ensemble.  We show that this quantity converges and describe its limiting distribution in terms of the Airy$_1$ random point field (i.e., the joint limit of the extremal eigenvalues of the GOE) from random matrix theory. 
}
\end{tabular}
\end{table}

\section{Introduction}

The 2-spin, spherical Sherrington-Kirkpatrick (SSK) model with zero magnetic field is defined by the random Hamiltonian
\begin{equation} \label{eqn:hamiltonian}
H_N(\sigma)=-\sum_{1\le i\neq j\le N} \frac{1}{\sqrt{2N}}g_{ij} \sigma_i\sigma_j.
\end{equation}
Here $H_N$ is a function of $\sigma\in \mathbb{S}^{N-1}:=\{\sigma\in \mathbb{R}^N, |\sigma|=\sqrt{N}\}$, and the coefficients $\{g_{ij}\}_{i, j=1}^N$ are iid standard normal random variables. 
 This model was introduced in \cite{KTJ}, by analogy with the standard SK model where the $\sigma$ are Ising spins taking values in the hypercube $\{\pm 1\}^N$ \cite{sherringtonkirkpatrick, panchenko}. The interest of spherical spin glass models lies in the availability of more explicit computations due to the continuous nature of the state space for the spins. See for example the papers \cite{subag1, subag2, subag3} by E. Subag for important work which takes advantage of the continuous geometry of spherical spin glasses.

The partition function of the SSK model is given by
\begin{equation}\label{eqn: partition}
Z_N(\beta)=\frac{1}{|\mathbb{S}^{N-1}|}\int_{\mathbb{S}^{N-1}}e^{-\beta H_N(\sigma)}\,\mathrm{d}\omega_N(\sigma).
\end{equation}
Here $\beta>0$ is a parameter corresponding to the inverse temperature. The $N\rightarrow \infty$ limit of the \emph{free energy}
\[F_N(\beta)=\frac{1}{N}\log Z_N(\beta)\]
was determined in \cite{CS}, and a rigorous justification appeared in \cite{talagrand-FE}. The model exhibits a phase transition at $\beta_c=1$, in the sense that the limit
\[F(\beta):= \lim_{N\rightarrow \infty} F_N(\beta)\]
fails to be analytic in $\beta$ at this value.

In \cite{baiklee}, J. Baik and J.O. Lee use a contour integral representation for the partition function which had previously appeared in \cite{KTJ} to compute the asymptotic fluctuations of $F_N$ around $F(\beta)$. 
They show that in the high temperature phase $\beta<\beta_c=1$, the quantity
\begin{equation}\label{eqn: gaussian}
  N(F_N(\beta)-F(\beta))
\end{equation}
converges to a normal random variable. This is the analogue for the SSK model of the classical central limit theorem for the SK model with Ising spins of M. Aizenman, D. Ruelle, and J. Lebowitz \cite{aizenmanlebowitzruelle}.   

In the low temperature phase of the SSK, Baik and Lee proved that the free energy has asymptotic fluctuations given by the Tracy-Widom distribution associated to the Gaussian Orthogonal Ensemble (GOE): 
\begin{equation}\label{eqn: low-temperature}
  \frac{N^{2/3}}{\beta-1}(F_N(\beta)-F(\beta))\Rightarrow TW_1.
\end{equation}
The convergence is in distribution as $N \to \infty$ and $TW_1$ denotes the asymptotic distribution of the top eigenvalue of the real symmetric GOE matrix \cite{tracywidom}. The analog of the low temperature result \eqref{eqn: low-temperature} for the classical SK model seems out of reach of current methods. Moreover, we are not aware of a prediction concerning the limiting distribution of the free energy of the SK model in the low temperature phase.

In a parallel development, a model related to \eqref{eqn:hamiltonian} and \eqref{eqn: partition} has appeared in the context of high-dimensional statistics.  Onatski, Moreira and Hallin \cite{stats_ssk} obtained an analog of the high-temperature CLT in the case that 
the random variables $\{g_{ij} \}_{ij}$ are associated with a Wishart ensemble (as opposed to the case under consideration where they are naturally associated with a symmetric matrix of normal random variables).  In this context, the Gaussian fluctuations have implications for the asymptotic power of statistical tests in detecting the presence of an unknown signal in an otherwise isotropically distributed dataset. Here, the high temperature regime corresponds to the regime of low signal-to-noise ratio.

In addition to the intrinsic interest of computing the fluctuations of $F_N(\beta)$, the method in \cite{baiklee} offers a satisfying interpretation of the phase transition in the SSK in terms of random matrix theory. The argument in \cite{baiklee} reveals that in the the high temperature phase $F_N(\beta)-F(\beta)$ is dominated by \emph{linear statistics} of eigenvalues $\{ \lambda_j \}_{j=1}^N$ of the matrix 
$$M_{ij}=\frac{1}{\sqrt{2N}}(g_{ij}+g_{ji}).$$
  Linear statistics are quantities of the form
\beq \label{eqn:xx1}
\sum_{j=1}^N f(\lambda_j)-N \int_{-2}^2 f(x) \rho_{\mathrm{sc}}(x)\,\mathrm{d}x,
\eeq
where $f$ is a regular function, and $\rhosc$ denotes the semicircle distribution, defined below.  This latter quantity is the asymptotic density of states of the GOE eigenvalues, and is known as Wigner's semicircle law \cite{wigner_orig}. The Gaussian behavior \eqref{eqn: gaussian} then follows from the well-known random matrix fact that the asymptotic fluctuations of \eqref{eqn:xx1} 
 are Gaussian \cite{pastur2011eigenvalue}. In the low temperature phase, $F_N(\beta)-F(\beta)$ instead depends to leading order only on the first eigenvalue $\lambda_1$. Thus the phase transition corresponds to a transition from a regime where all eigenvalues contribute to the limiting behavior, to one where only the leading eigenvalue does. Baik and Lee have applied their method to a number of variants of the SSK \cite{baiklee, baiklee1, baiklee2, baiklee3}, including the bipartite SSK and models incorporating a Curie-Weiss type interaction in addition to the spin glass couplings.

The phase transition in the classical and spherical SK model, and more generally in $p$-spin models \cite{panchenko}, can also be detected in the terms of the behavior of \emph{overlaps}. To define these, we introduce the Gibbs measure defined by the expectation
\begin{equation}\label{eqn: expectation}
  \langle f\rangle_{N,\beta} = \frac{1}{Z_N (\beta)} \int_{\mathbb{S}^{N-1}}e^{-\beta H_N(\sigma)} f(\sigma)\, \frac{\mathrm{d}\omega_N(\sigma) }{|\mathbb{S}^{N-1}|}
\end{equation}
Let $\sigma^{(1)},\sigma^{(2)}\in \mathbb{S}^{N-1}$ be two independent samples  (``\emph{replicas}'') from the Gibbs measure \eqref{eqn: expectation}. The overlap between $\sigma^{(1)}$ and $\sigma^{(2)}$ is the normalized inner product:
\[R_{12}=\frac{1}{N}(\sigma^{(1)}\cdot \sigma^{(2)}).\]
It is known that for $\beta<1$, $R_{12}$ tends to zero as $N\rightarrow \infty$, while in the low temperature phase $\beta>1$, it concentrates around the constant values $\pm q=\frac{1-\beta}{\beta}$ \cite{panchenkotalagrand}. In \cite{vl-ps}, V.-L. Nguyen and the second author used a contour integral representation related to that to that in \cite{baiklee1} to show that $R_{12}$ has Gaussian fluctuations for temperatures corresponding to $\beta\leq 1-N^{-1/3+\epsilon}$, for any $\epsilon >0$.

In the present work, we describe the annealed fluctuations of the Gibbs expectation $\langle R_{12}^2\rangle$ about the limiting value $\pm q$ in the low temperature regime. More precisely, Theorem \ref{thm1} below, provides an expansion for the overlap around these values down to order $o (N^{-2/3})$ in terms of explicit quantities related to the Gaussian Orthogonal Ensemble from random matrix theory. The expansion we derive does not seem to appear in the physics literature, but we were informed by J. Baik  that predictions close to the results we find were obtained using physics methods by Baik, Le Doussal and Wu \cite{baikprivate}.

The expansion of Theorem \ref{thm1} states that the leading order to the contribution of the fluctuations of the overlap around its mean is given by the quantity
\beq \label{eqn:xx2}
\frac{1}{N} \sum_{j=2}^N \frac{1}{ \lambda_j(M) - \lambda_1(M) }+1
\eeq
where the $\lambda_j (M)$ are the eigenvalues of $M$ arranged in decreasing order.  In our second main result, Theorem \ref{thm:mainconv}, we prove that (when renormalized by $N^{1/3}$) the quantity \eqref{eqn:xx2} converges and moreover describe its limit in terms of the Airy$_1$ random point field.  This latter point process arises in random matrix theory as the limits of the largest eigenvalues of the GOE.


\section{Main results}
We express the asymptotic distribution of the overlap between two replicas in terms of the eigenvalue distribution generated by the Gaussian Orthogonal Ensemble (GOE). To understand the connection between the GOE and our problem, define the symmetric random matrix $M$ defined by
\begin{align}M_{ij}=\begin{cases} 
\frac{g_{ij}+g_{ji}}{\sqrt{2 N}},& i\neq j\\
0,& i=j\\
\end{cases} \label{eqn:mdef}
\end{align}
where the $\{ g_{ij} \}_{i, j}$ are the random variables appearing in the definition of the Hamiltonian \eqref{eqn:hamiltonian}.  
The distribution of $M$ is that of a normalized GOE (Gaussian Orthogonal Ensemble) matrix with the diagonal set to zero. 
We denote the ordered eigenvalues of $M$ by
\[\lambda_1(M)\ge \lambda_2(M)\ge \ldots \ge \lambda_N(M).\]

Next, note that the Hamiltonian $H_N(\sigma)$ equals 
\begin{equation}\label{eqn: H-M}
\begin{split}
-\frac{1}{\sqrt{2N}}\sum_{1\le i\neq j\le N} g_{ij} \sigma_i\sigma_j&=-\frac{1}{2\sqrt{N}}\sum_{i\neq j}\frac{g_{ij}+g_{ji}}{\sqrt{2}} \sigma_i\sigma_j\\
&=: -\frac{1}{2}\langle \sigma, M\sigma\rangle.
\end{split}
\end{equation}
For two vectors $\sigma_1,\sigma_2\in \mathbb{S}^{N-1}$ we define the \emph{overlap} as the normalized inner product of $\sigma_1$ and $\sigma_2$:
\begin{align*}
R_{12}&=\frac{1}{N}(\sigma^{(1)}\cdot \sigma^{(2)})\\
&=\frac{1}{N}\sum_{i=1}^N \sigma_i^{(1)}\sigma_i^{(2)},
\end{align*}
where $(\sigma_i^{(k)})_{1\le i \le N}$, $k=1,2$ are the components of $\sigma^{(k)}$. For a bounded, measurable function 
\[f: (\mathbb{S}^{N-1})^k \rightarrow \mathbb{R},\]
we denote the Gibbs expectation of $f$ by
\[\langle f(\sigma^{(1)},\ldots,\sigma^{(k)})\rangle=\frac{1}{Z_N(\beta)^k}\frac{1}{|\mathbb{S}^{N-1}|^k}\int_{(\mathbb{S}^{N-1})^k}f(\sigma^{(1)},\ldots, \sigma^{(k)}) e^{-\beta \sum_{j=1}^k H_N(\sigma^{(j)})}\,\mathrm{d}\omega_N(\sigma^{(1)})\cdots \mathrm{d}\omega_N(\sigma^{(k)}).\]

As a consequence of the representation \eqref{eqn: H-M}, we derive the following integral formula for the Gibbs expectation  $\langle R_{12}^2 \rangle$ in Section \ref{sec: formulas}:
\beq\label{eqn: rep1}
\langle R_{12}^2 \rangle = \left(\frac{1}{4} \right) \frac{ \int \int \e^{N(Q(z) + Q(w) )/2} \left( \sum_{i=1}^N \frac{1}{ \beta^2 N^2 ( z-\lambda_i (M) ) ( w - \lambda_i (M) )}    \right) \d z \d w}{\left( \int \e^{ N Q(z)/2}  \d z \right)^2},
\eeq
where the integrals are over a vertical line  in the complex plane to the right of all the $\lambda_i (M)$ and 
\beq
Q(z) = \beta z - \frac{1}{N} \sum_i \log (z - \lambda_i (M)).
\eeq

Our main result provides an expansion of the overlap in terms of $\lambda_i(M)$, up to an error of size $o(N^{-2/3})$:
\begin{thm}\label{thm1} Let $\eps_1, \delta >0$ with $\frac{1}{3} > \delta$.  
For any $\eps >0$ and $N$ large enough, there is an event $\mathcal{F}_{\delta,\epsilon_1}$ such that
\[\mathbb{P}[\mathcal{F}_{\delta,\epsilon_1}]\ge 1-N^{-\delta+\epsilon}\]
on which the following estimate holds: 
\begin{align}
\langle R_{12}^2\rangle&=\left(\frac{1-\beta}{\beta}\right)^2+ 2\frac{\beta-1}{\beta^2} \cdot\left(\frac{1}{N} \sum_{j= 2}^N\frac{1}{\lambda_j (M)-\lambda_1(M)}+1\right) \nonumber \\
&-\frac{1}{N\beta^2}\left(\frac{1}{N}\sum_{j=2}^N \frac{1}{(\lambda_j(M)-\lambda_1(M))^2}\right)+\frac{1}{\beta^2}\left(\frac{1}{N} \sum_{j= 2}^N\frac{1}{\lambda_j(M)-\lambda_1(M)}+1\right)^2  \nonumber\\
&+\O(N^{3\delta+10\epsilon_1}N^{-1}). \label{eqn: rs}
\end{align}
\end{thm}
Theorem \ref{thm1} is a consequence of Theorem \ref{thm:main} below. The event $\F_{\delta, \eps_1}$ is defined in Definition \ref{def:Fevent} below; it is a high probability event on which certain a-priori estimates on the eigenvalue locations $\lambda_i$ hold (the rigidity and level repulsion estimates) - these are introduced in the next section.

Define
\beq
\tilm_N ( \lambda_1 ) = \frac{1}{N} \sum_{j\geq 2} \frac{1}{ \lambda_j (M)- \lambda_1(M)} , \qquad \tilm_N' ( \lambda_1 ) = \frac{1}{N} \sum_{j \geq 2} \frac{1}{ ( \lambda_j(M) - \lambda_1(M))^2}.
\eeq
Note that these quantities appear on the right side of equation \eqref{eqn: rs}.
The exponent $\delta>0$ is assocaited to \emph{level repulsion}, in that $\lambda_1 (M) - \lambda_2 (M) \geq N^{-2/3-\delta}$ on $\F_{\delta, \epsilon_1}$ (by definition of $\F_{\delta, \epsilon_1}$).   On the event $\F_{\delta, \epsilon_1}$  the magnitude of $\tilm_N$, $\tilm_N'$ will be seen to be at most
\beq
\tilm_N ( \lambda_1 ) + 1 = \O ( N^{-1/3+\delta+\eps_1} ), \qquad \frac{1}{N} \tilm_N' ( \lambda_1 ) = \O ( N^{-2/3+2\delta+\eps_1} ).
\eeq
Theorem \ref{thm1} thus identifies the overlap down to a term of order $N^{-2/3}$.   

It is interesting to study the behavior of the leading order contribution to the fluctuations of $\langle R_{12}^2 \rangle$ which is the term
\beq
\tilm_N ( \lambda_1 ) + 1,
\eeq
in the context of the work \cite{PT} of Panchenko and Talagrand.  They obtained an exponential estimate for the probability that $\langle R_{12}^2 \rangle \geq q^2 + \eps$ but noted that the event that $\langle R_{12}^2 \rangle \leq q^2 - \eps$ could not be ruled out at the level of large deviations.

Due to the rigidity estimates of random matrix theory which are reviewed in the next section, the quantity $\tilm_N ( \lambda_1 ) + 1$ has a light upper tail; for example the probability that it exceeds $N^{-1/3+\eps}$ goes to $0$ superpolynomially, for any fixed $\eps >0$.  On the other hand, the asymptotic density of the random variable $N^{2/3} ( \lambda_1 - \lambda_2)$ is expected to behave like $s^2$ near $0$ and so (recall that the eigenvalues are ordered and so this is a positive quantity) $\tilm_N ( \lambda_1) + 1$ has a relatively heavy lower tail.    Due to the somewhat large probability of the event $\F_{\delta, \eps_1}$ of Theorem \ref{thm1}, we do not attempt to make this comparison to the work \cite{PT} rigorous, settling for pointing out the heuristic agreement of our error term with the behavior observed in \cite{PT}.

As we have noted above, the matrix $M$ is closely related to the GOE.  The GOE is the matrix ensemble with entries
\begin{equation}\label{eqn: H-def}
\begin{split}
H_{ij}\sim \frac{1}{\sqrt{N}}N(0,1),&\quad 1\le i<j\le N,\\
H_{ii}\sim \frac{1}{\sqrt{N}}N(0,2),&\quad 1\le i \le N,\\
H_{ji}=H_{ij},&\quad 1\le i<j\le N,
\end{split}
\end{equation}
and all the non-identical random variables are independent.
We keep the dependence of $H$ on $N$ implicit.  In Appendix \ref{a:diag} we show that the largest  eigenvalues of $M$ agree with those of $H$ if we take $M_{ij} = H_{ij}$ for $i \neq j$, up to errors of order $\O(N^{-1+\eps})$ for any $\eps >0$.   Alternatively, it is possible to appeal to the literature on edge universality in random matrix theory (e.g., \cite{betaedge}), however the precise statement we require does not quite appear there.  Instead, we have opted to carry out the calculations in Appendix \ref{a:diag} which are a relatively straightforward application of the resolvent method.  The result derived here is in fact stronger than what could be deduced from the universality literature and may be of other application.

Related to this, we remark that the SSK Hamiltonian is for the most part defined as in \eqref{eqn:hamiltonian}, where we have excluded from the summation the diagonal terms $i = j$.  Of course in the usual SK model, whether or not the diagonal is included makes no difference since $\sigma_i^2=1$; in the SSK model, including the diagonal would result in $M = H$ above.  This would simplify our analysis somewhat as we could omit the calculations in Appendix \ref{a:diag} which compare the eigenvalues of $M$ directly to those of $H$.  To maintain consistency with the physics literature \cite{KTJ,crisanti} we have excluded the diagonal from the sum.

As can be seen by the expansion \eqref{eqn: rs}, the main contribution to the fluctuations of the overlap about its mean is from the extremal eigenvalues of $M$.  For any finite $k$, the largest $k$ eigenvalues of $H$, $\{ \lambda_i (H) \}_{i=1}^k$ are known to converge in distribution, after a rescaling, to the first $k$ particles of the Airy$_1$ random point field; we denote this latter quantity by $\{ \chi_i \}_{i=1}^{\infty}$.  We will give a more precise defintion in Section \ref{sec: convergence} below.  Due to our estimates proven in the appendix, the same joint convergence holds also for the largest eigenvalues of $M$.

A natural conjecture is then that the  rescaled fluctutations of the overlap converge in distribution to the random variable given by
\beq \label{eqn:xidef}
\Xi := \lim_{n \to \infty}  \left(  \sum_{j=2}^n \frac{1}{ \chi_j - \chi_1 } - \int_{0}^{ ( \frac{3 \pi n }{2} )^{2/3} } \frac{1}{ \pi \sqrt{ x} } d  x \right).
\eeq
In Theorem \ref{thm:conv1}, we show that this limit exists almost surely, and so $\Xi$ is a well-defined random variable.  The deterministic correction on the RHS of \eqref{eqn:xidef} represents the leading order term in the density of states of the Airy$_1$ random point field.  The expected location of the $j$th particle of the Airy$_1$ random point field  is roughly $\chi_j \sim j^{2/3}$ and so neither the sum or the deterministic correction converge as $n \to \infty$.  

Our main result on the limiting distribution of the fluctuations of the overlap is the following. 
\bet \label{thm:mainconv} Let $\Xi$ be the random variable in \eqref{eqn:xidef}.  
We have the following convergence in distribution for $\beta > 1$:
\beq
\lim_{N \to \infty} N^{1/3} \left[ \langle R_{12}^2 \rangle - \left( \frac{1-\beta}{\beta} \right)^2 \right]  = 2 \left( \frac{\beta-1}{\beta^2} \right) \Xi.
\eeq
\eet

In Theorem \ref{thm1} we  introduced the square in order to study the overlap, due to the symmetry of the overlap distribution with respect to the Gibbs measure (i.e., $\langle R_{12} \rangle = 0$).  An alternative would be to study $\langle | R_{12} | \rangle$; if we knew that $|R_{12}|$ concentrated about $q$ on the scale $N^{-1/3}$ then this of course could be deduced from Theorem \ref{thm1}.  We prove the concentration by calculating the fourth moment $\langle (R_{12}^2 - q^2 )^2 \rangle$; this is the content of the following theorem which is proven in Section \ref{sec:4m}.

\bet
On the event $\F_{\delta, \eps_1}$ of Theorem \ref{thm1} we have,
\beq \label{eqn:4mstatement}
\langle  ( R_{12}^2 - q^2 )^2 \rangle = \frac{ 8 ( \beta-1)^2}{ \beta^2} \frac{ \tilm_N'  ( \lambda_1 ) }{N} + 4 \frac{ ( \beta-1)^2}{\beta^4} ( 1 + \tilm_N ( \lambda_1 ) )^2 + \O (N^{-1+ 10 \eps_1+3 \delta}),
\eeq
and furthermore on the event $\F_{\delta, \eps_1}$, the first two terms are $\O ( N^{-2/3+2 \delta+10 \eps_1})$.  
As a consequence,
\beq
\langle |R_{12}| \rangle = q + \frac{1}{ \beta} ( \tilm_N ( \lambda_1) + 1 ) + \O ( N^{-2/3+2 \delta + 10\eps_1 })
\eeq
and so we have the convergence in distribution of
\beq
\lim_{N \to \infty} N^{1/3} \left[ \langle |R_{12} | \rangle - \frac{ 1 - \beta}{\beta} \right] = \frac{1}{ \beta} \Xi,
\eeq
where $\Xi$ is as above.
\eet

We discuss the relation of our results to the forthcoming work of Baik, Le Doussal and Wu \cite{baikprivate}.  They predict that the fluctuations of $R_{12}^2$ should be governed by
\beq
R_{12}^2 - q^2 \sim  \frac{ 2 ( \beta-1)}{\beta^2}  \frac{1}{N} \sum_{j=2}^N \frac{ n_j^2}{ \lambda_j - \lambda_1} + 1 
\eeq
where the $\{ n_j\}_j$ are independent standard normal random variables (in particular independent of the $\lambda_j$).  The Gibbs average $\langle \cdot \rangle $ corresponds to taking the expectation over the $\{ n_j \}_j$.  It is a simple calculation to integrate out the $\{n_j\}_j$ and find quantities agreeing with the  leading order contribution in \eqref{eqn: rs} and \eqref{eqn:4mstatement}.


\subsection{Outline of the paper}
In Section \ref{sec: rmt}, we state some basic results on random matrix theory which we will use to control the eigenvalues of the matrix $M$. In Section \ref{sec: formulas}, we obtain the representation \eqref{eqn: rep1}, along the lines of similar representations in \cite{baiklee1}, \cite{vl-ps}.

Most of the work is completed in Section \ref{sec: sd}, where we analyze the representation \ref{eqn: rep1} by the method of steepest descent. As was already noticed in \cite{baiklee1}, in the case $\beta>1$ of interest here, the analysis is complicated by proximity of the saddle point to the branch point $z=\lambda_1$ of the complex phase function $G(z)$ in \eqref{eqn: R12int}. Moreover, because the representation \eqref{eqn: rep1} involves a ratio of saddle point integrals, we must evaluate these to greater precision than was done in \cite{baiklee1} and the subsequent papers. The key idea is to separate the contribution from $\lambda_1$ to  $G(z)$ and use this to find the approximate steepest descent contours (see Lemma \ref{lem: saddle-lemma}).

In Section \ref{sec: convergence}, we show that the term $\tilm_N(\lambda_1)+1$ appearing in \eqref{eqn: rs} of order $N^{-1/3}$, in the sense that $N^{1/3}(\tilm_N(\lambda_1)+1)$ converges in distribution. This involves some establishing some preliminary estimates for the GOE as well as the Airy$_1$ random point field, which we could not locate in previous literature.  We deduce this using corresponding results for the GUE and Airy$_2$ random point field proven by Gustavsson and Soshnikov \cite{soshnikov,gust}, respectively, and the  Forrester-Rains coupling \cite{FR} between the GUE and GOE.

\section{Random matrix results}\label{sec: rmt}
In this section, we summarize the results from random matrix theory we use in the rest of the paper. A central role is played by the \emph{resolvent} matrix
\beq
R(H, z) = \frac{1}{H-z},  \qquad R(M, z) = \frac{1}{M-z}
\eeq
where $H$ is the GOE matrix in \eqref{eqn: H-def}, and $M$ is the matrix ensemble given by \eqref{eqn:mdef}.   The spectral parameter is $z$ is commonly denoted $z=E+\i \eta$ with $E, \eta \in \rr$ and $\eta>0$. In the recent literature, the resolvent has customarily been denoted by $G$, a notation we reserve for the quantity \eqref{eqn: G-def} in this paper. We also introduce the Stieltjes transform of the empirical eigenvalue distribution:
\[m_N(H, z)=\frac{1}{N}\mathrm{tr}\,R(H, z)=\frac{1}{N}\sum_{j=1}^N\frac{1}{\lambda_j(H)-z},\]
and similarly for $M$.  The classical semi-circle law is then equivalent to the approximation for fixed $z$,
\[m_N(H, z)=m_{\mathrm{sc}}(z)+o(1),\]
where the semi-circle law and its Stieltjes transform are
\[\rhosc (E) = \frac{1}{ 2 \pi  } \sqrt{ (4 - E^2)_+} , \qquad m_{\mathrm{sc}}(z)=\int \frac{1}{x-z}\rho_{\mathrm{sc}}(x)\,\mathrm{d}x.\]
We now state the local semi-circle law as it appears in \cite[Theorem 2.6]{KBG}. First, we introduce the notion of \emph{overwhelming probability}.
\bed
We say that an event or family of events $\{ \A_i \}_{i \in I}$ hold with overwhelming probability if for all $D>0$ we have $\sup_{i \in I} \pp [ A_i^c] \le N^{-D}$ for $N$ large enough.
\eed
\begin{thm}[Local semi-circle law]\label{thm: sc-law}
Define the spectral domain $\mathbf{S}$ by
\[\mathbf{S}=\{E+i\eta: |E|\le 10, 0<\eta\le 10\}.\]
For any $\epsilon>0$ and all $N$ sufficiently large, we have for both $R(z) = R(M, z)$ and $R(H, z)$ that the estimates,
\[\max_{i,j}|R_{ij}(z)-\delta_{ij}\msc(z)|\le \sqrt{\frac{\Im \msc (z)}{N^{1-\epsilon}\eta}}+\frac{1}{N^{1-\epsilon}\eta},\]
and for both $m_N(z) = m_N (H, z)$ and $m_N (M, z)$,
\[|m_N(z)-m_{\mathrm{sc}}(z)|\le \frac{1}{N^{1-\epsilon}\eta}\]
 hold uniformly in $z\in \mathbf{S}$ with overwhelming probability.
\end{thm}
A consequence of the semi-circle we will use several times is that the eigenvalues $\lambda_i$ are close to the corresponding quantiles of the semi-circle distribution. These quantiles are known as the \emph{classical locations} of the eigenvalues in random matrix theory:
\begin{equation}\label{eqn: gammai-def}
\int^{\gamma_i}_{-2} \rho_{\mathrm{sc}}(x)\,\mathrm{d}x=\frac{i}{N}.
\end{equation}
\begin{thm}[Eigenvalue rigidity] \label{thm: rigidity} For each $\epsilon>0$, we have that the estimates
\[|\lambda_i-\gamma_i|\le N^{-2/3+\epsilon}\min\{i, (N+1-i)\}^{-1/3}\]
hold uniformly in $i$ with overwhelming probability, for $\lambda_i$ the eigenvalues of $M$ or $H$.
\end{thm}
We will also need some finer information concerning \emph{level repulsion}. The next result shows that, up to an $\O(N^\epsilon)$ error, the distribution of the spacing between $\lambda_1$ and $\lambda_2$ has a density on scale $N^{-2/3}$. While this will be sufficient for our purposes, one instead expects that there is \emph{level repulsion}, i.e. $s$ on the right side of \eqref{eqn: s} should be replaced by $s^2$. This has been established in great generality for the spacings $\lambda_j-\lambda_{j+1}$ where $j\gg 1$ in \cite[Theorem 3.7]{edgebeta}, but has not been proven for the eigenvalues at the edge.   

The following result could be deduced from Remark 1.5 of \cite{ky}.  A complete proof was not given in that work and relies on asymptotics of the Hermite polynomials.  For the sake of completeness, we will give a different proof which relies only on the eigenvalue rigidity and the loop equations.  
\begin{lem}[Existence of spacing density] \label{lem:lr}
Let $\epsilon>0$. There is a constant $C>0$ such that for $N^{-1/3+\epsilon} \le s\le 1$,
\begin{equation}\label{eqn: s}
\mathbb{P}(N^{2/3}(\lambda_1-\lambda_2)<s)\le CsN^\epsilon,
\end{equation}
where the $\lambda_i$ are the eigenvalues of $M$ or $H$.  
\end{lem}
\remark Inspecting the proof we see that the restriction $s \geq N^{-1/3+\epsilon}$ enters only in proving the estimate for $M$ - i.e., it holds for all $s$ for the eigenvalues of $H$.

\begin{proof}
By Proposition \ref{prop:ev} it suffices to prove the estimate for $\lambda_i := \lambda_i (H)$.  
We begin with the obvious estimate:
\[\frac{1}{\lambda_1-\lambda_2}\le \sum_{j=2}^{N^\epsilon} \frac{1}{\lambda_1-\lambda_j}.\]
As a consequence of Theorem \ref{thm: rigidity} (see, for example \cite[Eqn (6.3)]{baiklee}), we have with overwhelming probability,
\[0\le \frac{1}{N}\sum_{j=2}^N\frac{1}{\lambda_1-\lambda_j}-\frac{1}{N}\sum_{j=2}^{N^\epsilon} \frac{1}{\lambda_1-\lambda_j}=1-\O(N^{-1/3+\epsilon}).\]
Combining this with the estimate $\ee[ | \lambda_2 - \lambda_1 |^{-3/2} ] \leq N^C$ for some $C >0$, which is a consequence of Section 5 of \cite{landonyau} we obtain the same inequality in expectation,
\[0\le \ee  \left[\frac{1}{N}\sum_{j=1}^N\frac{1}{\lambda_1-\lambda_j}-\frac{1}{N}\sum_{j=2}^{N^\epsilon} \frac{1}{\lambda_1-\lambda_j} \right]=1-\O(N^{-1/3+\epsilon}).\]
Then using Markov's inequality we have,
\begin{align*}\mathbb{P}(N^{2/3}(\lambda_1-\lambda_2)<s)&\le \mathbb{P}\left(\sum_{j=2}^{N^\epsilon} \frac{1}{\lambda_1-\lambda_j}>\frac{N^{2/3}}{s}\right)\\
&\le sN^{1/3}\mathbb{E}\big[ \frac{1}{N}\sum_{j=1}^{N^\epsilon} \frac{1}{\lambda_1-\lambda_j}\big]\\
&\le sN^{1/3} \left(\mathbb{E}[\frac{1}{N}\sum_{j=1}^{N} \frac{1}{\lambda_1-\lambda_j}-1\big]+\O(N^{-1/3+\epsilon})\right).
\end{align*}
By \cite[Lemma 3.7]{GS}, we have
\[\mathbb{E}[\frac{1}{N}\sum_{j=1}^{N} \frac{1}{\lambda_1-\lambda_j}\big]=\frac{\mathbb{E}[\lambda_1]}{2} = 1 + \O (N^{-2/3+\epsilon} )\]
where in the last step we used again Theorem \ref{thm: rigidity}.  The result follows.
\end{proof}\qed

\section{Representation for the overlap}\label{sec: formulas}
In this section, we derive a contour integral representation for Gibbs expectation $\langle R_{12}^2 \rangle$.  Throughout this section, we will denote the eigenvalues of $M$ by 
\beq
\lambda_i := \lambda_i (M),
\eeq
for notational simplicity. 
We now prove the following lemma. 
\begin{lem}\label{lem: representation} The quantity 
$\langle R_{12}^2\rangle$ is given by
\begin{equation}\label{eqn: R12int}
\langle R_{12}^2 \rangle = \frac{\int_{\gamma-\i\infty}^{\gamma+\i\infty} \int_{\gamma-\i\infty}^{\gamma+\i\infty}e^{\frac{N}{2}(G(z)+G(w))}\bigg(\sum_{i=1}^N \frac{1}{\beta^2N^2(z-\lambda_i)(w-\lambda_i)}\bigg)\mathrm{d}z \mathrm{d}w}{\Big (\int_{\gamma-\i\infty}^{\gamma+\i\infty} e^{\frac{N}{2}G(z)} dz\Big )^2},
\end{equation}
where 
\begin{equation}\label{eqn: G-def}
G(z)= \beta z - \frac{1}{N}\sum_{i=1}^N \log (z-\lambda_i),
\end{equation}
for any $\gamma \in \rr$ so that $\gamma > \lambda_1$. 
\end{lem}
\remark Note that up to constants the quantity appearing in the denominator of \eqref{eqn: R12int} is the partition function and is the representation used by Baik and Lee \cite{baiklee}.

\begin{proof}
Our starting point is the definition, 
\begin{equation}
\langle R_{12}^2\rangle= \frac{1}{Z_N(\beta)^2} \int_{( \mathbb{S}^{N-1})^2} \exp\bigg(\frac{\beta}{2} \langle \sigma^{(1)}, M\sigma^{(1)}\rangle+\frac{\beta}{2} \langle \sigma^{(2)} , M\sigma^{(2)} \rangle  \bigg)\big(\frac{1}{N}\langle \sigma^{(1)} ,\sigma^{(2)} \rangle\big)^2 \mathrm{d}\omega_N(\sigma^{(1)}) \mathrm{d} \omega_N(\sigma^{(2)}).
\end{equation}
Baik and Lee give the representation:
\begin{equation}\label{eqn: partition-rep}
Z_N(\beta)=\frac{\Gamma(N/2)\cdot 2^{N/2-1}}{2\pi i (N\beta)^{N/2-1}}\int_{\gamma-i\infty}^{\gamma+i\infty} e^{\frac{N}{2}G(z)}\,\mathrm{d}z.
\end{equation}
Let $S^{N-1} = \{ \mathbf{x} \in \mathbb{R}^N: \|\textbf{x}\| =1\}$ be the unit sphere in $\mathbb{R}^N$. Let $\mathrm{d}\Omega$ be the surface area measure on $S^{N-1}$, so that $ \vert S^{N-1} \vert^{-1}\mathrm{d}\Omega$ is the uniform measure on $S^{N-1}$. By a change of variables, we obtain
\begin{align}
&\int_{( \mathbb{S}^{N-1})^2} \exp\bigg(\frac{\beta}{2} \langle \sigma^{(1)}, M\sigma^{(1)} \rangle+\frac{\beta}{2} \langle \sigma^{(2)}, M\sigma^{(2)}\rangle  \bigg) (\frac{1}{N}\langle \sigma^{(1)},\sigma^{(2)} \rangle)^2\mathrm{d}\omega_N(\sigma^{(1)}) \mathrm{d} \omega_N(\sigma^{(2)})\\
=&\frac{1}{| S^{N-1}|^2}\int_{( S^{N-1})^2} \exp\bigg(\frac{\beta}{2} N \langle \mathbf{x}, M\mathbf{x}\rangle+\frac{\beta}{2} N\langle \mathbf{y}, M\mathbf{y}\rangle \bigg) \langle \mathbf{x},\mathbf{y} \rangle^2 \mathrm{d}\Omega_1 \mathrm{d} \Omega_2.
\end{align}
Let $z,w\in  \{u \in \cc: \Re\, u>\lambda_1(M)\}$. In order to compute the above integral, we consider 
\begin{align*}\label{def:J}
  J(z,w)&=\int_{\mathbb{R}^N}\int_{\mathbb{R}^N} e^{\frac{\beta}{2} N\sum_{i=1}^N (\lambda_i-z) x_i^2} e^{\frac{\beta}{2} N\sum_{i=1}^N (\lambda_i -w) y_i^2} \big(\sum_{i=1}^N x_iy_i\big)^2\,\prod_{i=1}^N \mathrm{d}x_i\mathrm{d}y_i\\
  &=\int_{\mathbb{R}^N}\int_{\mathbb{R}^N} e^{\frac{\beta}{2} N\sum_{i=1}^N (\lambda_i-z) x_i^2} e^{\frac{\beta}{2} N\sum_{i=1}^N (\lambda_i -w) y_i^2} \sum_{i=1}^N x_i^2y_i^2\,\prod_{i=1}^N \mathrm{d}x_i\mathrm{d}y_i.
\end{align*}
We use polar coordinates, substituting $\textbf{x}=r_1\textbf{x}_1$ and $\mathbf{y}=s_1\mathbf{y}_1$ with $r_1,s_1>0$ and $\|\mathbf{x}_1\|=\|\mathbf{y}_1\|=1$. We then set $(\beta/2) Nr_1^2=r$, $(\beta/2) Ns_1^2=s$
to find that 
\begin{align*}
 J(z,w)=\frac{2^N}{(\beta N)^N}\int_0^\infty \int_0^\infty e^{- z r}e^{- w s} I(r,s)s^{\frac{N}{2}-1}r^{\frac{N}{2}-1}\,\mathrm{d}r\mathrm{d}s,
\end{align*}
where
\[ I(r,s)= \int_{S^{N-1}\times S^{N-1}} e^{ r\langle \mathbf{x}_1, M\mathbf{x}_1\rangle +s\langle \mathbf{y}_1, M\mathbf{y}_1\rangle}   \frac{rs}{\beta^2N^2}\langle\mathbf{x}_1, \mathbf{y}_1 \rangle ^2 \,\mathrm{d}\Omega_1 \mathrm{d}\Omega_2.\]
On the other hand, direct integration shows that the function $J$ is given by
\begin{align}
J(z,w)= \big (\frac{2\pi}{\beta N}\big )^N \prod_{i=1}^N \frac{1}{\sqrt{(z-\lambda_i)(w-\lambda_i)}} \bigg(\sum_{i=1}^N \frac{1}{\beta^2N^2(z-\lambda_i)(w-\lambda_i)}\bigg).
\end{align}
Taking the inverse Laplace transform, we obtain
\begin{align}
&\frac{2^N}{(\beta N)^N} I(r,s)s^{\frac{N}{2}-1}r^{\frac{N}{2}-1}\\
=&\frac{1}{(2\pi i)^2}\int_{\gamma-i\infty}^{\gamma+i\infty} \int_{\gamma-i\infty}^{\gamma+i\infty} e^{zr}e^{ws} J(z,w) \mathrm{d}z \mathrm{d}w\\
=&\Big (\frac{2\pi}{\beta N}\Big )^N \frac{1}{(2\pi i)^2}\int_{\gamma-i\infty}^{\gamma+i\infty} \int_{\gamma-i\infty}^{\gamma+i\infty} e^{zr}e^{ws} \prod_{i=1}^N \frac{1}{\sqrt{(z-\lambda_i)(w-\lambda_i)}} \bigg(\sum_{i=1}^N \frac{1}{\beta^2N^2(z-\lambda_i)(w-\lambda_i)}\bigg) \mathrm{d}z \mathrm{d}w
\end{align}
where $\gamma$ is any real number satisfying $\gamma >\lambda_1$. Recalling that
\[|S^{N-1}|=\frac{2\pi^{N/2}}{\Gamma(\frac{N}{2})}\]
and letting $r=s=\frac{\beta}{2} N$, we obtain:
\begin{align}
&\frac{1}{|S^{N-1}|^2}\int_{(S^{N-1})^2} e^{ N\frac{\beta}{2}\langle \mathbf{x}_1, M\mathbf{x}_1 \rangle +N\frac{\beta}{2}\langle \mathbf{y}_1, M\mathbf{y}_1\rangle}   \langle\mathbf{x}_1,\mathbf{y}_1\rangle^2\,\mathrm{d}\Omega_1 \mathrm{d}\Omega_2 \nonumber \\=&\frac{2^{N-2}\Gamma(N/2)^2}{(2\pi i)^2 (\beta N)^{N-2}}\int_{\gamma-i\infty}^{\gamma+i\infty} \int_{\gamma-i\infty}^{\gamma+i\infty} e^{N\frac{\beta}{2}(z+w)}\prod_{i=1}^N \frac{1}{\sqrt{(z-\lambda_i)(w-\lambda_i)}} \bigg(\sum_{i=1}^N \frac{1}{\beta^2N^2(z-\lambda_i)(w-\lambda_i)}\bigg) \mathrm{d}z \mathrm{d}w. \label{eqn: final-rep}
\end{align}
Combining \eqref{eqn: final-rep} and \eqref{eqn: partition-rep}, we obtain \eqref{eqn: R12int}. \qed
\end{proof}

\section{Steepest descent analysis}\label{sec: sd}
We proceed to the asymptotic evaluation of the integrals in \eqref{eqn: R12int}.   As in the previous section, we will continue to denote the eigenvalues of $M$ by
\beq
\lambda_i := \lambda_i (M).
\eeq
As was already noticed in \cite{baiklee}, in the low temperature regime, the dominant contribution to the integrals comes from an $O(N^{-1})$  neighborhood of the saddle point $\gamma$ which is itself  distance $N^{-1}$ from the largest eigenvalue of $M$.  The prescence of the branch point due to $\lambda_1$ close to the saddle makes a steepest descent analysis via a direct expansion of the function $G(z)$ untenable.

Compared to the computation in \cite{baiklee} and subsequent works, we must  evaluate the numerator and denominator in \eqref{eqn: R12int} with greater precision. For the main result of \cite{baiklee}, for example, it was sufficient to show that
\[ Z(\beta) = iKe^{\frac{N}{2}G(\gamma)},\]
where $K$ satisfies $N^{-C} \leq K \leq N^C$ for any $C>0$.  The contribution of $K$ to the free energy is then $\O (N^{-1} \log(N) )$, automatically of lower order than the dominant Tracy-Widom fluctuations which are of size $\O (N^{-2/3})$. 

In order to evaluate the overlap, it is necessary to determine the leading order term of $Z ( \beta)$, not only up to multiplicative terms. Additionally, our computation involves a more precise localization of $\gamma$ than $|\gamma -\lambda_1|\le N^{-1+\epsilon}$.

We now give an overview of the saddle point analysis.    The saddle point of the function $G(z)$ is distance of order $N^{-1}$ from $\lambda_1.$  Instead of working directly with the steepest descent contours of the function $G(z)$, we will consider the dominant contribution near the saddle which is, up to additive constants,
\begin{equation}\label{eqn: simplified}
(\beta-1) z -\frac{1}{N}\log(z-\lambda_1).
\end{equation}
This function is much simpler than $G(z)$, as it involves only the eigenvalue $\lambda_1$.   The contributions from other eigenvalues are replaced by their deterministic leading order term using Theorem \ref{thm: rigidity}.  The additional key input here is Lemma \ref{lem:lr} which ensures that the eigenvalues $\{ \lambda_j \}_{j=2}^N$ are an order of magnitude further from $\lambda_1$ than the distance between the saddle and $\lambda_1$.  This allows for the localization of the function $G(z)$ near its saddle, despite the prescence of the branch point due to the logarithmic singularity at $ z= \lambda_1$.

  The saddle point of the function \eqref{eqn: simplified} is clearly,
\beq \label{eqn: saddle1}
\gamma := \lambda_1 + \frac{\cb}{N}
\eeq
where 
\beq \label{eqn: saddle2}
\cb := \frac{1}{\beta-1}.
\eeq
The advantage afforded by working with the approximation \eqref{eqn: simplified} is that the behavior of the steepest descent contours of this function are relatively explicit.  
For the contours, we make the ansatz $z=\gamma + E+i\eta(E)$, for $E \leq 0$.  Setting the imaginary part of \eqref{eqn: simplified} to zero gives a parametrization of the approximate steepest descent contour we will use.  We determine properties of this parameterization in Lemma \ref{lem: saddle-lemma}.  In Lemma \ref{lem:reg} we analyze the behavior of $G(z)$ along this approximate steepest descent contour.  
\bel \label{lem: saddle-lemma}
For $z\in \mathbb{C}\setminus (-\infty,0]$, let $-\pi < \mathrm{arg}(z)\le \pi$ be the standard determination of the argument.
For $E <0$ the equation
\beq \label{eqn:eta}
\eta ( \beta-1) = \frac{1}{N} \,\mathrm{arg} (E + \i \eta + \cb/N )
\eeq
has a unique strictly positive solution which we denote $\eta(E)$.   Furthermore, there is a constant $c_1>0$ so that if $0 \geq E \geq - c_1/N$,
\beq \label{eqn:etaest1}
N \eta(E) = \sqrt{ 3 \cb  |NE| } ( 1 + \O ( |N E| ) ).
\eeq
For any $c_2>0$ there is a $c_3 >0$ depending on $c_2>0$ so that if $E \leq - \frac{c_2}{N}$, we have 
\beq \label{eqn:etaest2}
\frac{c_3}{N} \leq \eta \leq  \frac{\pi}{N (\beta-1) }.
\eeq

\eel
Before proceeding to the proof of Lemma \ref{lem: saddle-lemma}, we note that if $f$ is analytic with real and imaginary parts denoted by
\beq
f (x + \i y ) = u (x, y) + \i v (x, y),
\eeq
then, with $\del_z = (\del_x - \i \del_y)/2$, the Cauchy-Riemann equations imply
\beq
f'(z) = u_x - \i u_y = v_y + \i v_x
\eeq
and so
\beq\label{eqn: CR}
\del_E \Re [f] = \Re[ f'], \quad \del_{\eta} \Re[f] = - \Im [f'], \quad \del_E \Im[f] = \Im[f'], \quad \del_\eta \Im[f] = \Re[f'],
\eeq
using the notation $z = E + \i \eta$.

\proof For uniqueness we note that if $E < - \cb/N$, then the left side of \eqref{eqn:eta} is increasing, whereas the right side is decreasing.  If $E > - \cb/N$ we calculate the derivative of the right side,
\beq
\del_\eta \frac{1}{N} \arg (E + \i \eta + \cb/N) = \frac{1}{N} \frac{ E + \cb/N}{(E+\cb/N)^2 + \eta^2 }.
\eeq
This is a decreasing function of $\eta$, and so the right side of \eqref{eqn:eta} is a concave function of $\eta$. Its derivative at $\eta =0$ is strictly greater than $(\beta-1)$, so we get the uniqueness as the left side is a linear function with slope $(\beta-1)$.  Differentiating the equation \eqref{eqn:eta}, we find using \eqref{eqn: CR}
\beq
\frac{ \d \eta }{\d E} \left( \beta -1 - \frac{1}{N} \frac{E+\cb /N}{(E+\cb /N)^2+ \eta^2 } \right) = \frac{1}{N} \frac{ - \eta}{(E+\cb /N)^2 + \eta^2} .
\eeq
Note that the second factor on the left is positive (for $E> - \cb /N$ it is the difference of the slopes of the tangent lines of the functions on either side of \eqref{eqn:eta} at the point $\eta(E)$), so 
\beq
\frac{ \d \eta}{\d E} \leq 0,
\eeq
and the lower bound of \eqref{eqn:etaest2} will follow once we establish \eqref{eqn:etaest1}.  The upper bound is immediate.

Let $\log$ denote the principal determination of the logarithm.  Expanding this function in a power series around $\cb/N$, we have for some $c>0$ and $|z| \leq c/N$,
\beq \label{eqn:expa1}
\log (z + \cb/N) = z \frac{N}{\cb} - \frac{z^2}{2} \frac{N^2}{\cb^2}+ \frac{z^3}{6} \frac{2 N^3}{\cb^3} + N^4 z^4 f (z)
\eeq
where $f(z)$ is an analytic function in the disc $|z| \leq c/N$, obeying the estimates
\beq
| \Im[ f(z) ]|  \leq CN  \Im [z] , \qquad |f(z) | \leq C.
\eeq
These estimates follow from the fact that all the coefficients in the power series expansion of the logarithm are real.
The imaginary part of $\log$ is the argument function appearing on the right side of \eqref{eqn:eta}.   Taking imaginary parts on both sides of \eqref{eqn:expa1} and using \eqref{eqn:eta}, we find (denoting $\eta = \eta (E)$ for brevity)
\beq \label{eqn: solvethis}
0 = -E \eta \frac{N^2}{\cb^2} + E^2 \eta\frac{N^3}{\cb^3} - \frac{\eta^3}{3}\frac{N^3}{\cb^3} + \O \left( N^4 \eta E^3 + N^4 E \eta^3 \right)
\eeq
Dividing by $N\eta$ gives,
\beq
N \eta \leq C \sqrt{ N |E|}
\eeq
after possibly making $c>0$ smaller.  Using this to estimate the higher order terms in \eqref{eqn: solvethis} yields \eqref{eqn:etaest1} after solving \eqref{eqn: solvethis} for $N \eta$ as a function of $NE$. 
\qed

We also require the following elementary lemma.
\bel \label{lem:rega1}
Let $\eta (E)$ be as above.  For any $c_1 >0$, there is a $c_2 >0$ so that if $ E \leq - c_1 / N$ then,
\beq \label{eqn:rega1}
\beta -1 - \frac{1}{N} \frac{ E + \cb/N}{ (E+ \cb/N)^2 + \eta(E)^2 } \geq c_2.
\eeq
\eel
\proof The case $E+ \cb/N \leq 0$ is trival so we may assume $E \geq - \cb/N$.  Recall that $\eta (E)$ is the unique positive solution to $(N \eta ) / \cb = \arctan (\eta/ (E + \cb/N ) )$.  If $NE/\cb + 1  \leq \frac{1}{10}$, then $\arctan ( 1 / ( NE/\cb+1)) \geq \pi/3 > 1$ and so for such $E$, we have $N \eta (E) / \cb \geq 1$.  We have,
\beq
\frac{1}{ \cb} - \frac{1}{N} \frac{ E + \cb/N}{(E+ \cb/N)^2+( \eta(E))^2 } = \frac{ (N \eta(E))^2 + ( NE + \cb)(NE) }{ \cb ( (N E + \cb )^2 + ( N \eta (E) )^2 )}
\eeq
The denominator is bounded above, and the numerator is bounded below by $\cb^2/2$ for $NE / \cb + 1 \leq \frac{1}{10}$.  It remains to consider the case where $NE/ \cb +1 \in (c_3, 1-c_3 )$ for fixed $c_3 >0$.  Consider the function,
\beq
f ( \eta ) =  \frac{1}{N}\arctan ( \eta / (E + \cb/N)  - \eta / \cb.
\eeq
This function has zeros at $\eta = 0$ and $\eta = \eta (E)$ and is strictly positive in between these points.  We need to prove that there is a constant $c'$ depending on $c_3$ so that $f' ( \eta (E) )  < - c'$.   By direct calculation, it has a local maximum in between these two points at $N \eta^* = \sqrt{ (NE + \cb)(-NE)}.$  We see that $f' ( 0) > c$ for a $c>0$ and that $0 \geq f'' ( \eta) \geq - CN$ for $\eta \geq 0$ and constants $c, C$ depending on $c_3$.  Since $f' ( \eta^* ) = 0$ it follows that $\eta^* \geq c /N$ for some new $c>0$ depending on $c_3$, and then that $f ( \eta^* ) \geq c'/N$ for some $c' >0$.  Since $f' ( \eta)$ is bounded, we then see that $\eta (E) - \eta^* \geq c''/N$ for some $c''>0$.   Since $\eta^* \geq c/N$ we see that $f''(\eta) \leq - N c''$ for $\eta(E) \geq \eta \geq \eta^*$ for some $c''' >0$.  This then implies that $f' ( \eta(E)) \leq - c_4$ for some $c_4 >0$, which is what we needed to prove.
\qed

We now define the event $\F_{\delta, \eps_1}$ of Theorem \ref{thm1}.

\begin{definition} \label{def:Fevent}
Let $\frac{1}{3} > \delta >0$ and $\eps_1 >0$.  Let $\F_{\delta, \eps_1}$ be the following event:
\beq
\F_{\delta, \eps_1} = \left\{ N^{2/3} (\lambda_1 - \lambda_2 ) > N^{-\delta} \right\} \cap \left\{ |\lambda_i - \gamma_i | \leq \frac{ N^{\eps_1/10}}{ \min\left\{ i^{1/3}, (N+1-i)^{1/3} \right\} N^{2/3} }, i=1,\ldots,N \right\},
\eeq
where $\gamma_i$ are the classical eigenvalue locations defined in \eqref{eqn: gammai-def}.  
\end{definition}
As stated above, this is the event in the statement of Theorem \ref{thm1}.  By Theorem \ref{thm: rigidity} and Lemma \ref{lem:lr} we have that
\beq
\pp [ F_{\delta, \eps_1} ] \geq 1 - N^{-\delta+\eps'}
\eeq
for any $\eps' >0$ and $N$ large enough.  Fix a sufficiently small $\kappa >0$.  In particular, we take 
\beq
\kappa < \frac{1/3-\delta}{10}
\eeq
Define the contours
\beq
\Gamma_1 := \{ E \pm \i \eta(E) : 0 \geq E \geq - N^{-1+\kappa } \}
\eeq
and
\beq
\Gamma_2 := \{ -N^{-1+\kappa} \pm \i \eta : \eta \geq \eta (-N^{-1+\kappa} ) \}
\eeq
The contour $\Gamma_1$ is a U-shaped contour symmetric along the negative real axis, and $\Gamma_2$ are vertical lines at the ends of $\Gamma_1$ going to $ \pm \i \infty$.
\bel \label{lem:reg}
Assume that $\F_{\delta, \eps_1}$ holds with $\eps_1$ sufficiently small.  The following estimates hold.  There is a $c>0$ so that 
\beq \label{eqn:reg2}
\Re [ G (z + \gamma) ] - \Re[ G (\gamma) ] \leq -c N^{-1+\kappa}, \quad z = N^{-1+\kappa} + \i \eta,\quad \eta\geq 0.
\eeq
For any $c_1 >0$ there are $c_2 >0$ and $C_1$ so that the following holds for $z \in \Gamma_1$ and large enough $N$:
\begin{align} 
\Re [ G (z + \gamma) ] - \Re[ G (\gamma) ]  &\leq \1_{\{ |E| \leq c_1/N \}} |E| N^{-1/3+\eps_1+\delta} \notag\\
-&\1_{\{ |E| \geq c_1/N \} } \left[ (|E|-c_1/N )c_2 -  C_1 N^{-1-1/3+\eps_1+\delta}\right], \label{eqn:reg1}
\end{align}
Finally, there is the following estimate for $-N^{-1+\kappa} \leq E \leq 0$ and $\eta \geq 10$,
\beq \label{eqn:reg3}
\Re[G(z+\gamma)]  - \Re[G (z) ] \leq -\frac{1}{3} \log(1+ \eta).
\eeq
\eel
\proof We calculate some derivatives of $\Re\,G$ along the contours.  In this proof, $z$ will be restricted to lie on the various contours and so we will generally denote $z = E + \i \eta (E)$.    First, along $\Gamma_1$,
\begin{align}
\frac{\d}{\d E} \Re\, G ( \gamma + E + i \eta(E) ) &= \beta + \Re[ m_N(\gamma+z) ] - \Im [ m_N(\gamma+z) ] \frac{ \d \eta }{ \d E} \notag\\
&\geq \beta + \Re[m_N(\gamma+z)],
\end{align}
where we used that $\d \eta / \d E$ is negative.  We write,
\beq
\Re [m_N(\gamma+z)] = - \frac{1}{N} \frac{ \cb/N + E}{(\cb/N +E)^2 + \eta^2 } - \frac{1}{N} \sum_{j=2}^N \frac{E + \gamma - \lambda_j }{ ( E + \gamma - \lambda_j )^2 + \eta^2 }.
\eeq
We need to estimate the second term.  We write,
\beq \label{eqn: split-real}
\frac{1}{N} \sum_{j=2}^N \frac{E + \gamma - \lambda_j }{ ( E + \gamma - \lambda_j )^2 + \eta^2 }  = \frac{1}{N} \sum_{j=2}^{\lfloor N^{\eps_1}\rfloor }  \frac{E + \gamma - \lambda_j }{ ( E + \gamma - \lambda_j )^2 + \eta^2 }  + \frac{1}{N} \sum_{j=\lfloor N^{\eps_1}\rfloor +1}^N  \frac{E + \gamma - \lambda_j }{ ( E + \gamma - \lambda_j )^2 + \eta^2 } .
\eeq
From the level repulsion assumption and choice of $\kappa$,
\beq
|E + \gamma - \lambda_j | \geq N^{-2/3-\delta} - N^{-1+\kappa} \geq c N^{-2/3-\delta},
\eeq
for all $j \geq 2$, and so
\beq
\left| \frac{1}{N} \sum_{j=2}^{N^{\eps_1}}  \frac{E + \gamma - \lambda_j }{ ( E + \gamma - \lambda_j )^2 + \eta^2 } \right| \leq C N^{\eps_1+\delta-1/3}.
\eeq
For the second term in \eqref{eqn: split-real}, rigidity gives
\beq
\left| \frac{1}{N} \sum_{j=\lfloor N^{\eps_1}\rfloor+1}^N  \frac{E + \gamma - \lambda_j }{ ( E + \gamma - \lambda_j )^2 + \eta^2 }  - \Re\left[  \int_{\gamma_{N^{\eps_1}}}^{2} \frac{ \rho_\mathrm{sc} (x) \d x }{ (\gamma+z)-x}  \right] \right| \leq C N^{-1/3+\eps_1/10}.
\eeq
Now, since $|\gamma+z - 2| \leq C N^{\eps_1/10-2/3}$,
\beq
\left|  \int_{\gamma_{N^{\eps_1}}}^{2} \frac{ \rho_\mathrm{sc}}{ x - (\gamma+z)} -  \int_{\gamma_{N^{\eps_1}}}^{2} \frac{ \rho_\mathrm{sc}}{ x - 2} \right| \leq N^{\eps_1/10-2/3} C \int_{N^{-2/3}}^{2} \frac{\sqrt{x}}{x^2} \leq C N^{\eps_1-1/3}.
\eeq
Finally,
\beq
\left| \int^{2}_{\gamma_{N^\eps_1}} \frac{\rho_{\mathrm{sc}}(x)}{x-2} \d x -1\right| \leq N^{\eps_1-1/3}.
\eeq
Therefore,
\beq
\beta + \Re[m_N(\gamma+z)] \geq \left( \beta-1 -\frac{1}{N} \frac{\cb/N+E}{(\cb/N+E)^2 + \eta^2 } \right) - C N^{\eps_1+\delta-1/3}.
\eeq
As observed in the proof of Lemma \ref{lem: saddle-lemma}, the first term on the RHS (in the brackets) is positive.  If $E \leq -c_1/N$ for a $c_1>0$ then by Lemma \ref{lem:rega1} we conclude that  there is a $c_2 >0$ depending on $c_1$ so that 
\beq
\beta + \Re[m_N(\gamma+z)]  \geq c_2.
\eeq
We have therefore proven that for $-c_1/N \leq E \leq 0$,
\beq
\frac{ \mathrm{d} }{ \mathrm{d} E } \Re[ G ( E + \i \eta (E) + \cb/N)] \geq - C N^{\eps_1+\delta-1/3}
\eeq
and for $-N^{-1+\kappa} \leq E \leq - c_1/N$,
\beq
\frac{ \mathrm{d} }{ \mathrm{d} E } \Re[ G ( E + \i \eta (E) + \cb/N)]  > c_2/2.
\eeq
The estimate \eqref{eqn:reg1} follows from the previous two estimates and integration.  We consider $z$ of the form $z = -N^{-1+\kappa} + \i \eta$, in order to prove \eqref{eqn:reg2}.  We  consider the behavior of $\Re[G (z) ]$ as $\eta$ varies.  We calculate, 
\beq \label{eqn:ff1}
\del_\eta \Re[G(\gamma + z)] = - \Im [ m_N ( \gamma+z)].
\eeq
This is decreasing, so we immediately get -- using \eqref{eqn:reg1} -- the estimate \eqref{eqn:reg2} in the region $\eta \geq \eta (-N^{-1+\kappa})$.  For smaller $\eta \leq \eta ( - N^{-1+\kappa} )$, note that 
\beq
| \gamma -N^{-1+\kappa}-\lambda_1 | \geq c N^{-1+\kappa}, \quad | \gamma -N^{-1+\kappa}-\lambda_2 | \geq cN^{-2/3-\delta} \geq c N^{-1+\kappa}.
\eeq
Hence,
\beq
\Im [ m_N] \leq N^{\eps_1} \frac{N \eta}{N^{2 \kappa}}+ \frac{1}{N} \sum_{j \geq N^{\eps_1}} \frac{ \eta}{|\gamma + z - \lambda_j|^2} \leq N^{\eps_1} \frac{N \eta}{N^{2 \kappa}}+\eta N^{1/3}.
\eeq
Since $\eta (-N^{-1+\kappa} ) \leq C/N$, we get \eqref{eqn:reg2} for the rest of the possible values of $\eta$, integrating $\Re[ G]$ from $z = - N^{-1+\kappa} + \i \eta (N^{-1+\kappa} )$ to a smaller $\eta \leq \eta (N^{-1+\kappa} )$, using \eqref{eqn:ff1} and the above estimate on the derivative.  

Finally, we turn to \eqref{eqn:reg3}.  We have,
\beq
\Re[ G(\gamma + z ) ] - \Re [ G ( \gamma) ] = \beta E -\frac{1}{N} \sum_j \log \left| 1 + \frac{z}{ \gamma-\lambda_j } \right|.
\eeq
For $j \leq N^{\eps_1}$,
\beq
\left| 1 + \frac{z}{ \gamma-\lambda_j } \right| \geq \eta N^{1/3}
\eeq
and for the rest of $j$,
\beq
\left| 1 + \frac{z}{ \gamma-\lambda_j } \right|  \geq \frac{\eta}{3}
\eeq
This yields \eqref{eqn:reg3} (recall $\eta \geq 10$ for this estimate). \qed

Now introduce the contour
\beq
\hat{\Gamma} = \Gamma_1 \cup \Gamma_3,
\eeq
where
\beq
\Gamma_3 = \{ -N^{-1+\kappa} + \i \eta : 0  \leq\eta \leq \eta (-N^{-1+\kappa } ) \}.
\eeq
Recall that $\Gamma_1$ is a U-shaped contour symmetric about the negative real axis.  The contour $\Gamma_3$ connects the ends of the U to the real axis. 

The following lemma shows that we can deform the contours in \eqref{eqn: R12int} into $\hat{\Gamma}$.
\bel \label{lem: hatcontours}
Suppose the event $\F_{\delta, \eps_1}$ holds with $\eps_1$ sufficiently small.  Then the following estimates hold.  There is a $c>0$ so that
\begin{align} \label{eqn:ct1}
 & \int_{\gamma-\i \infty}^{\gamma+ \i \infty} \int_{\gamma-\i \infty}^{\gamma+ \i \infty} \e^{\frac{N}{2} (G(z) + G(w) - 2 G ( \gamma)) } \left( \sum_{i=1}^N \frac{1}{ \beta^2 ( z +\gamma- \lambda_i ) (w +\gamma- \lambda_i ) } \right) \notag\\
=  & \int_{\hat{\Gamma} \times \hat{\Gamma}}\e^{\frac{N}{2} (G(z) + G(w)  - 2 G( \gamma)) } \left( \sum_{i=1}^N \frac{1}{ \beta^2 ( z +\gamma- \lambda_i ) (w+\gamma - \lambda_i ) } \right)  + \O ( \e^{ - N^c } )
\end{align}
and
\begin{align} \label{eqn:ct2}
\int_{\gamma-\i \infty}^{\gamma+ \i \infty} e^{ \frac{N}{2} (G(z)- G( \gamma ) ) } = \int_{\hat{\Gamma}} \e^{\frac{N}{2} ( G (z) - G ( \gamma ))} + \O ( \e^{ - N^c} ).
\end{align}
\eel
\proof By analyticity and the absolute convergence guaranteed by \eqref{eqn:reg3}, all of the contours can be moved from the vertical lines appearing in \eqref{eqn: R12int}  to the contour $\Gamma_1 \cup \Gamma_2$.  It just remains to replace $\Gamma_2$ by $\Gamma_3$.    This replacement for \eqref{eqn:ct2} is immediate from \eqref{eqn:reg2} and \eqref{eqn:reg3}.  For \eqref{eqn:ct1} we also have to deal with cross terms (e.g, the integral over $\Gamma_2$ or $\Gamma_3$ in the $z$ variable times the integral over $\Gamma_1$ in the $w$ variable) and the extra terms in the integrand.  Note that along all of the contours under consideration we always have
\beq
|z  + \gamma - \lambda_j | \geq N^{-1}.
\eeq
Note furthermore that for $z \in \Gamma_1$, the estimate \eqref{eqn:reg1} gives
\beq
N ( \Re G (z) - \Re G (\gamma ) ) \leq C,
\eeq
and that the arc length of $\Gamma_1$ satisfies
\beq
| \Gamma_1 | \leq C N^{-1+\kappa}
\eeq
as $\eta(E)$ is monotonic.  These observations together with \eqref{eqn:reg3} and \eqref{eqn:reg2} yield \eqref{eqn:ct1}. \qed

In the following lemma we Taylor expand $G(z)$ (or at least all of the terms appearing in its definition except the one with $\lambda_1$) around the saddle $\gamma$, arriving at a form of the integrands which we will be able to calculate.  We first define a few functions which naturally appear in the Taylor expansion.  Let,
\beq
\tilm_N (w) := \frac{1}{N} \sum_{j=2}^N \frac{1}{ \lambda_j - w }
\eeq
and
\beq
g(z) := (\beta + \tilm_N ( \gamma ) ) z - \frac{1}{N} \log(1+Nz / \cb )
\eeq
\bel \label{lem: hattaylor}
The following holds on the event $\F_{\delta, \eps_1}$ for sufficiently small $\eps_1>0$.   First, we have
\begin{align} \label{eqn:bb3}
\int_{\hat{\Gamma}} \e^{\frac{N}{2} (G(z+\gamma) - G ( \gamma ))} = \int_{\hat{\Gamma}} \e^{\frac{N}{2} g(z) }(1  + N z^2 \frac{\tilm_N' ( \gamma)}{4}  ) + \O \left( N^{-2+4 \kappa+3\delta+\eps_1} \right)
\end{align}
Second,
\begin{align} 
&\int_{\hat{\Gamma} \times \hat{\Gamma}} \e^{\frac{N}{2}  (G(z + \gamma) + G(w+\gamma) - 2 G ( \gamma ))} \frac{1}{ N^2 ( z+\cb/N )(w+\cb/N) } \notag\\
= &\int_{\hat{\Gamma} \times \hat{\Gamma} } \e^{ \frac{N}{2} ( g(z) + g ( w)) } \left(1 + N z^2 \frac{ \tilm_N' ( \gamma)}{4} \right) \left(1 + N w^2 \frac{ \tilm_N' ( \gamma)}{4} \right) \frac{1}{ N^2 ( z+\cb/N )(w+\cb/N) } \notag\\
+ & \O\left( \frac{ N^{5 \kappa+\eps_1+3 \delta}}{N^3} \right). \label{eqn:bb1}
\end{align}
Finally,
\begin{align}
&\int_{\hat{\Gamma} \times \hat{\Gamma}} \e^{\frac{N}{2}  (G(z + \gamma) + G(w+\gamma) - 2 G ( \gamma ))} \sum_{j=2}^N \frac{1}{ N^2 ( z+\gamma -\lambda_j )(w+\gamma-\lambda_j) } \notag\\
=& \frac{1}{N^2} \sum_{j=2}^N \frac{1}{ (\lambda_1-\lambda_j )^2} \left( \int_{ \hat{\Gamma}} \e^{ \frac{N}{2} g(z) } \left(1 + N z^2 \frac{ \tilm_N' ( \gamma ) }{4} \right) \right)^2+  \O\left( \frac{ N^{5 \kappa+\eps_1+3 \delta}}{N^3} \right). \label{eqn:bb2}
\end{align}
\eel
\proof We write
\begin{align}
G(z+\gamma) - G ( \gamma ) &= \beta z -\frac{1}{N} \log \left( 1 + \frac{ N z}{ \cb} \right) \notag\\
&- \frac{1}{N} \sum_{j\geq 2} \log(z + \gamma - \lambda_j )  - \log ( \gamma- \lambda_j )
\end{align}
We Taylor expand the second term.  For $z \in \hat{\Gamma}$,
\begin{align}
\frac{1}{N} \sum_{j \geq 2} \log (z + \gamma - \lambda_j ) - \log ( \gamma - \lambda_j ) = - z \tilm_N ( \gamma ) - \frac{z^2}{2} \tilm_N' ( \gamma ) + \O\left(\frac{N^{3 \kappa+\eps_1 + 3 \delta}}{N^2} \right).
\end{align}
From the fact that $z^2 \tilm_N' ( \gamma ) = \O ( N^{-5/3+5\kappa+2 \delta+ \eps_1 } )$ for $z \in \hat{\Gamma}$ we first conclude that
\begin{align}
\Re [g (z) ] \leq \frac{C}{N},
\end{align}
where we used Lemma \ref{lem:reg} (i.e, the corresponding estimate for $\Re[G (z + \gamma) ] - \Re[ G ( \gamma ) ]$ and the above two equalities which relate this quantity to $g(z)$).  Then, for $z \in \hat{\Gamma}$,
\begin{align} \label{eqn:aa1}
\e^{\frac{N}{2} (G(z+\gamma) - G ( \gamma ) )} = \e^{ \frac{N}{2} g(z) } \left(1 +N  \frac{z^2}{4} \tilm'_N ( \gamma ) \right) +  \O\left(\frac{N^{3 \kappa+\eps_1 + 3 \delta}}{N} \right).
\end{align}
Note that
\begin{align} \label{eqn:aa2}
\left| \frac{1}{N^2} \sum_{j=1}^N \frac{1}{(z-\lambda_j )(w - \lambda_j ) } \right| \leq C,
\end{align}
and
\begin{align} \label{eqn:aa3}
\left| \frac{1}{N^2} \sum_{j=2}^N \frac{1}{ (\lambda_j - z-\gamma) (\lambda_j - w-\gamma) } - \frac{1}{N^2} \sum_{j=2}^N \frac{1}{ ( \lambda_j - \lambda_1)^2 } \right| \leq C \frac{N^{\kappa+3 \delta+\eps_1}}{N}.
\end{align}
The first and second estimates of the lemma follows from \eqref{eqn:aa1} and the fact that $| \hat{\Gamma} | \leq C N^{-1+\kappa}$.   For the third estimate, one first uses  \eqref{eqn:aa2} and \eqref{eqn:aa1} to arrive at an integral in terms of $g(z), g(w)$ and the (quantity inside the absolute value) on the LHS of \eqref{eqn:aa2}.  The error is $\O (N^{5 \kappa + \eps_1+3 \delta-3} )$.  The next replacement uses \eqref{eqn:aa3}, and one arrives at the final estimate of the lemma. \qed

We now rescale and shift the contour of integration to lie along the real axis.  Let $\Gamma_{r}$ be the following keyhole contour around the point $ -\cb$, for $r < \cb/10$:
\beq
\Gamma_r := \{ E \pm \i 0 : E < - \cb - r \} \cup \{ z : |z- \cb| = r \}.
\eeq
\bel \label{lem: contours-final}
On the event $\F_{\delta, \eps_1}$ we have, for sufficiently small $\eps_1 >0$, the following estimates for some $c>0$.
\begin{align}
\int_{\hat{\Gamma}} \e^{ \frac{N}{2} g(z) } \left( 1 + N z^2 \frac{ \tilm_N' ( \gamma) }{4} \right) = \frac{1}{N} \int_{\Gamma_r} \frac{ \e^{( \beta + \tilm_N ( \gamma ) ) u/2 } }{ \sqrt{1+u/ \cb } } \left( 1 +\frac{u^2 \tilm_N' ( \gamma ) }{ 4 N }  \right) \d u+ \O ( \e^{ - N^c} )
\end{align}
and
\begin{align} \label{eqn:2mf2}
\int_{\hat{\Gamma}} \e^{ \frac{N}{2} g(z) } \left(1 + N z^2 \frac{ \tilm_N' ( \gamma ) }{4} \right) \frac{1}{ N ( z + \cb/N )} &= \frac{1}{N} \int_{\Gamma_r} \frac{\e^{ ( \beta+\tilm_N ( \gamma )) u/2}}{ \sqrt{1+u/\cb}}\left( 1 + \frac{ u^2\tilm_N' ( \gamma )}{4N} \right) \frac{ \d u}{ u + \cb } \notag\\
&+ \O ( \e^{ - N^c} )
\end{align}

\eel
\proof First we make the substitution $u = z/N$.  As the integrand is analytic on $\mathbb{C} \backslash \{ E  \leq - \cb \}$ we see that all of the contours may be shifted from $N \hat{\Gamma}$ to $ \Gamma_r \backslash \{ E \pm \i 0 :  E \leq - N^{\kappa } \}$.  The rest of $\Gamma_r$ may be added to the integral as $\beta + \tilm_N ( \gamma ) > c_1$ for some $c_1 >0$, at only an error exponentially small in $N$. \qed

We collect some explicit integrals in the next lemma.
\bel \label{lem:explicitintegrals}
Let $a, b>0$, and $\Gamma_{r, b}$ be a keyhole contour around $-b$ as above.  Then,
\begin{align}
\int_{\Gamma_{r, b}}  \frac{ \e^{ a z}}{ \sqrt{ z + b} } &= 2 \i \frac{ \e^{-ab}}{\sqrt{a}} \sqrt{\pi}\\
\int_{\Gamma_{r, b}} \e^{az} \sqrt{z+b}  &= - \i \e^{-ab} a^{-3/2} \sqrt{\pi} \\
\int_{\Gamma_{r, b} }\e^{ a z } (z+b)^{3/2} &= \frac{3}{2} \i  \e^{-ab} a^{-5/2}  \sqrt{ \pi }\\
\int_{\Gamma_{r, b}}  \frac{ \e^{ a z}}{ \sqrt{ z + b} } z^2 &= \frac{ \sqrt{\pi} \e^{-ab} \i }{\sqrt{a}} \left( \frac{3}{2 a^2}+ \frac{2b}{a}+ 2 b^2 \right)  \\
\int_{\Gamma_{r, b}}  \frac{ \e^{ a z}}{( z + b)^{3/2}}  &= \sqrt{a} \e^{-ab} 4 \i \sqrt{\pi} \\
\int_{\Gamma_{r, b}}  \frac{ \e^{ a z}}{ (z+b)^{3/2}}  z^2&= \frac{ \sqrt{\pi} \e^{-ab} \i}{ \sqrt{a}} \left( \frac{-1}{a}-4b+ 4b^2 a \right)
\end{align}
\eel
All of the above calculations can be done by considering the contributions from the integral along the real axis and the circle around $z=b$ as $r \to 0$. In the cases where these contributions are diverging, one treats the circular integral by Taylor expansion (i.e., expanding the exponential around $z=b$), and integrates by parts the integral along the real axis.  One finds that the diverging quantities cancel, and is left with a real integral which can be calculated explicitly.

We finally arrive at the following, from which Theorem \ref{thm1} follows. 
\bet \label{thm:main} On the event $\F_{\delta, \eps_1}$ we have with sufficiently small $\eps_1$ that,
\begin{align}
\langle R_{12}^2 \rangle =\left( 1 - \beta^{-1} \right)^2 + 2 \frac{ \beta-1}{\beta^2} \left( \tilm_N ( \lambda_1) + 1 \right) -  \frac{ \tilm_N' ( \lambda_1 ) }{N \beta^2} + \frac{ ( \tilm_N ( \lambda_1 ) +1 )^2}{ \beta^2} + \O \left( \frac{ N^{ 3 \delta+ 10 \eps_1 }}{N} \right).
\end{align}
\eet

\begin{proof}  We first use Lemma \ref{lem: representation} to arrive at the formula \eqref{eqn: R12int} for the overlap.  The results of the present section are used to analyze the contour integrals appearing in the numerator and denominator of \eqref{eqn: R12int}.  From \eqref{eqn:ct1}, \eqref{eqn:bb1}, \eqref{eqn:bb2} and Lemma \ref{lem: contours-final} we arrive at the following expression for the numerator:
\begin{align}
& \int_{\gamma-i\infty}^{\gamma+i\infty} \int_{\gamma-i\infty}^{\gamma+i\infty}e^{\frac{N}{2}(G(z)+G(w))}\bigg(\sum_{i=1}^N \frac{1}{\beta^2N^2(z-\lambda_i)(w-\lambda_i)}\bigg)\mathrm{d}z \mathrm{d}w \notag\\
= & \left(  \frac{1}{N \beta} \int_{\Gamma_r} \frac{\e^{ ( \beta+\tilm_N ( \gamma )) u/2}}{ \sqrt{1+u/\cb}}\left( 1 + \frac{ u^2\tilm_N' ( \gamma )}{4N} \right) \frac{ \d u}{ u + \cb }  \right)^2 \notag\\
+ & \frac{1}{(N\beta)^2} \sum_{j=2}^N \frac{1}{ ( \lambda_1 - \lambda_j )^2} \left( \frac{1}{N} \int_{\Gamma_r} \frac{ \e^{( \beta + \tilm_N ( \gamma ) ) u/2 } }{ \sqrt{1+u/ \cb } } \left( 1 +\frac{u^2 \tilm_N' ( \gamma ) }{ 4 N }  \right) \d u \right)^2 + \O ( N^{5 \kappa + \eps_1 + 3 \delta-3} ) \label{eqn:r12final1}
\end{align}
 for $\kappa$ as above.  For the integral in the denominator we use \eqref{eqn:ct2}, \eqref{eqn:bb3} and Lemma \ref{lem: contours-final} to find,
 \begin{align}
 \int_{\gamma-i\infty}^{\gamma+i\infty} e^{\frac{N}{2}G(z)} dz = \frac{1}{N} \int_{\Gamma_r} \frac{ \e^{( \beta + \tilm_N ( \gamma ) ) u/2 } }{ \sqrt{1+u/ \cb } } \left( 1 +\frac{u^2 \tilm_N' ( \gamma ) }{ 4 N }  \right) \d u + \O (N^{3 \kappa+\eps_1-2} ).
 \end{align}
 We now use Lemma \ref{lem:explicitintegrals} with $a = \frac{1}{2} ( \beta+ \tilm_N ( \gamma) )$ and $b = \cb$.  For the numerator, the two terms of \eqref{eqn:r12final1} equal
 \begin{align}
 \left( \frac{2}{ \sqrt{a}} \frac{ \sqrt{ \cb} \e^{ - a \cb } \i \sqrt{ \pi}}{N \beta} \right)^2 \left( (2 a )^2 - \frac{ \tilm_N' ( \gamma) }{ N} + \frac{ \tilm_N' ( \lambda_1) }{N} \right) + \O ( N^{5 \kappa + \eps_1 + 3 \delta-3} )
 \end{align}
 whereas the denominator equals
 \begin{align}
  \left( \frac{2}{ \sqrt{a}} \frac{ \sqrt{ \cb} \e^{ - a \cb } \i \sqrt{ \pi}}{N} \right)^2  \left(1 + \frac{ 3 \tilm_N' ( \gamma)}{4 N a^2 } \right)
 \end{align}
 We get the claim from these two calculations as well as,
 \beq \label{eqn:r12final2}
 \tilm_N ( \gamma ) = \frac{1}{N} \sum_{j=2}^N \frac{1}{ \lambda_j - \lambda_1 - \cb/N} = \frac{1}{N} \sum_{j=2}^N \frac{1}{ \lambda_j - \lambda_1 } + \frac{\cb}{N^2} \sum_{j=2}^N \frac{1}{ ( \lambda_j - \lambda_1 )^2 } + \O \left( N^{\eps_1+3 \delta-1} \right).
 \eeq
 \qed
\end{proof}

\section{Existence of limit}
\label{sec: convergence}

In this section we consider the limit of the random variables
\beq
-N^{1/3}(\tilde{m}_N(\lambda_1)+1)=N^{1/3} \left( \frac{1}{N} \sum_{j=2}^N \frac{1}{ \lambda_1 - \lambda_j } - 1 \right)
\eeq
as $N \to \infty$.  The limit will be characterized in terms of the Airy$_1$ random point field.  Our convention is so that if $\lambda_1 \geq \lambda_2 \geq \dots$ are the largest eigenvalues of the GOE, then for every finite $k$,
\beq
\{ N^{2/3} ( 2 - \lambda_j ) \}_{j=1}^k \to \{ \chi_j \}_{j=1}^k,
\eeq
so that the ensemble $\{ \chi_j \}_{j=1}^\infty$ has finitely many particles located on the negative real line.  We will prove the following theorem.  
\bet \label{thm:conv1}
Let $\lambda_1 \geq \lambda_2 \geq \dots$ denote the largest eigenvalues of the GOE, and $\{ \chi_j \}_{j=1}^\infty$ the Airy$_1$ random point field.  The sequences of random variables
\beq
N^{1/3} \left( \frac{1}{N} \sum_{j=2}^N \frac{1}{ \lambda_1 - \lambda_j } - 1 \right)
\eeq
converges in distribution to a random variable $\Xi$ which is given by
\beq \label{eqn:airya1}
\Xi = \lim_{n \to \infty} \left( \sum_{j=2}^n \frac{1}{ \chi_j - \chi_1 } - \frac{1}{\pi} \int_0^{( \frac{3 \pi n}{2} )^{2/3}} \frac{d x}{\sqrt{x}} \right),
\eeq
where the limit on the RHS of \eqref{eqn:airya1} exists almost surely.
\eet
\remark  We do not determine whether or not the distribution of $\Xi$ is non-trivial.

\subsection{Preliminary estimates}

We will need an estimate for the variance of the number of eigenvalues of the GOE in an interval as well as the corresponding estimate for the Airy$_1$ random point field.  We will deduce these from the corresponding results for the GUE and Airy$_2$ random point field and the coupling of Forrester and Rains between the GUE and the GOE \cite{FR}.  
\bet [Soshnikov \cite{soshnikov}] Let $\chit_i$ be the particles of the Airy$_2$ point process and let $T>0$.  We have the following estimates for some $C>0$ and any $T>0$.
\beq
\left| \ee \left[ \left| \{ i : \chit_i \leq T \} \right| \right] - \frac{2}{3 \pi} T^{3/2} \right|  \leq C,
\eeq
and
\beq \label{eqn: soshnikov-var}
\left| \mathrm{Var} \left(  \left| \{ i : \chit_i \leq T \} \right|\right) - \frac{3}{ 4 \pi^2} \log T \right| \leq C.
\eeq
\eet

\remark Soshnikov states the variance asymptotics \eqref{eqn: soshnikov-var} with the constant $\frac{11}{12\pi^2}$ instead of $\frac{3}{4\pi^2}$. This appears to be due to a mistake in the computation of the quanity $I_3(u)$ in \cite[Lemma 5]{soshnikov}. Note also that the factor $3/4$ is consistent with the variance asymptotics for the counting function in the GUE, in \eqref{eqn:vargue} below.  We will use only that the variance grows logarithmically in $T$.

\bet [Gustavsson \cite{gust}] Let $\{ \mu_i \}_{i=1}^N$ be the eigenvalues of the GUE.  Let $\eps >0$.  There is a $C>0$ so that the following holds.  For any $0 \leq s \leq N^{2/3-\eps}$, we have
\beq \label{eqn:expgue2}
\left| \ee\left[\left| \{ i : \mu_i \geq 2 - s N^{-2/3} \} \right| \right] -N \int_{2- s N^{-2/3}}^2 \frac{ \sqrt{4 - E^2} \d E}{2 \pi} \right| \leq C
\eeq
and so for $s \leq N^{4/15}$, 
\beq \label{eqn:expgue}
\left| \ee\left[\left| \{ i : \mu_i \geq 2 - s N^{-2/3} \} \right| \right] - \frac{2}{3 \pi} s^{3/2} \right| \leq C.
\eeq
Furthermore, there is a $C_1 >0$ so that if $C_1 \leq s \leq N^{2/3-\eps}$, then
\beq \label{eqn:vargue}
\left| \mathrm{Var} \left( \left| \{ i : \mu_i \geq 2 - sN^{-2/3} \} \right| \right) - \frac{3}{4  \pi^2} \log (s) \right| \leq C ( 1 +  \log \log (s) )
\eeq
\eet
\remark Gustavsson only claimed the result \eqref{eqn:vargue} in the case that $s \to \infty$ as $N \to \infty$ at any arbitrarily slow rate (his interest was in the case that the variance tends to infinity, a necessary condition for applying a theorem of Costin and Lebowitz \cite{costinlebowitz}).  Inspecting his proof yields the estimate \eqref{eqn:vargue} for fixed but large enough $s$. \qed

We need also the following result of Forrester and Rains.
\bet [Forrester, Rains, \cite{FR}]  \label{thm:fr} Let $\mathrm{GOE}_n$ and $\mathrm{GUE}_{n}$ denote the set formed by the union of the eigenvalues of the GOE and GUE, respectively.  Then,
\beq 
\mathrm{GUE}_n \stackrel{d}{=} \mathrm{Even} \left( \mathrm{GOE}_n \cup \mathrm{GOE}_{n+1} \right)
\eeq
where the RHS is the set formed by the second largest, fourth largest, sixth largest, etc. elements of $\mathrm{GOE}_n \cup \mathrm{GOE}_{n+1} $.
\eet

From the above results we deduce the following.
\bep
Let $\{ \lambda_i \}_{i=1}^N$ denote the eigenvalues of the GOE and let $\eps >0$.  There is a $C>0$ so that the following holds.  For any $0 \leq s \leq N^{4/15}$, we have
\beq \label{eqn:expgoe}
\left| \ee \left[ \left| \{ i : \lambda_i \geq 2 - s N^{-2/3} \} \right| \right] - \frac{2}{3 \pi} s^{3/2} \right| \leq C.
\eeq
There is a $C_1 >0$ so that if $C_1 \leq s \leq N^{2/3-\eps}$, then
\beq \label{eqn:vargoe}
\mathrm{Var} \left( \left| \{ i : \lambda_i \geq 2 - s N^{-2/3} \} \right|  \right) \leq C \log(s)
\eeq
\eep
\proof Let $\{ \lambda_i^{(N)} \}_{i=1}^N$ and $\{ \lambda_i^{(N+1)} \}_{i=1}^{N+1}$ be the eigenvalues of two independent GOE matrices of dimension $N$ and $N+1$ respectively.   For the course of this proof, let $X^{(1)}_N (s)$ and  $X^{(1)}_{N+1} (s)$ be the number of eigenvalues of these matrices at least $2 - sN^{-2/3}$.  Let $X^{(2)}_N (s)$ be same but for the eigenvalues of an independent GUE matrix.  The coupling of Theorem \ref{thm:fr} implies that there is a random variable $Y$ and a bounded random variable $Z$ so that,
\beq
X^{(2)}_N(s) \stackrel{d}{=} Y, \qquad Y- Z = \frac{1}{2} \left( X^{(1)}_N (s) + X^{(1)}_{N+1} (s) \right) .
\eeq
Using that $X^{(1)}_N (s)$ and $X^{(1)}_{N+1} (s)$ are independent, that $Z$ is bounded, and the estimate \eqref{eqn:vargue} yields \eqref{eqn:vargoe}.  Taking expectations we see that
\beq \label{eqn:expcoup}
\ee[ X_N^{(2)}(s) ] = \frac{1}{2} \left( \ee[ X_N^{(1)} (s)] + \ee[ X_{N+1}^{(1)} (s) ] \right) + \O (1).
\eeq
We now estimate the difference of the two quantities on the RHS.  Given a GOE matrix $H$ of dimension $N+1$, its minor formed by removing the first row and column is a GOE matrix of dimension $N$ multiplied by the prefactor $a_N = \sqrt{N/(N+1)}$.   Additionally, the $N$ eigenvalues of the minor interlace the $N+1$ eigenvalues of $H$.  Hence,
\beq
\left| \ee[ X_{N+1}^{(1)}(s) ] - \ee[ X_{N}^{(1)} ( s + (2  - s N^{-2/3} )N^{2/3}(1-a_N^{-1} )  ]\right| \leq C.
\eeq
Let $b_N =  -N^{2/3}(2 -  s N^{-2/3} )(1- a_N^{-1} )= \O ( N^{-1/3} )$.  Applying now \eqref{eqn:expcoup} twice, with $s$ and $s-b_N$, and taking the difference we see that (note that the difference $X^{(1)}_n (s) - X_{n}^{(1)}(s-b_N)$ with $n=N, N+1$ has the same sign as it is just the number of eigenvalues between $s$ and $s-b_N$) 
\beq
\left| \ee[ X_{N}^{(1)} (s)] - \ee[ X_N^{(1)} (s- b_N ) ] \right| \leq \left| \ee[ X_{N}^{(2)} (s)] - \ee[ X_N^{(2)} (s- b_N ) ]  \right| + C \leq C
\eeq
where we used \eqref{eqn:expgue} in the last step.  This yields \eqref{eqn:expgoe}. \qed

\bel
Let $\{ \chi_i \}_{i}$ and $\{ \chi'_i \}_i$ be two independent Airy$_1$ random point fields.  Let $\{ \zeta_j \}_j$ be an Airy$_2$ random point field.  Let $T>0$.  Then,
\beq
\pp\left[ | \{ i : \zeta_i \leq T \} | = k \right] 
\eeq
is equal to the probability that there are either $2k$ or $2k+1$ particles from the superimposed point process $\{ \chi_i \}_i \cup \{ \chi'_i \}_i$ below $T$.
\eel
\proof Let $\mu_i$ be the scaled eigenvalues of a GUE matrix, and $\lambda_i$ and $\lambda_i'$ the scaled eigenvalues from independent GOE matrices of dimension $N$ and $N+1$ (for the latter, do the scaling with $N^{-2/3}$ and not by $(N+1)^{-2/3}$).  By the convergence in distribution of the $k$th eigenvalue of the GOE/GUE to the $k$th particle of the Airy$_1$/Airy$_2$ process we have,
\begin{align}
\pp \left[ | \{ i : \zeta_i \leq T \} | = k \right]  &= \pp[ \zeta_k \leq T, \zeta_{k+1} > T ] = \pp [ \zeta_k \leq T ] - \pp [ \zeta_{k+1} \leq T ] \notag\\
=& \lim_{N \to \infty} \pp [ \mu_k \leq T ] - \pp [ \mu_{k+1} \leq T ] = \lim_{N \to \infty} \pp [ \left| \{ i : \mu_i \leq T \} \right|= k ].
\end{align}
By Theorem \ref{thm:fr},
\begin{align}
 \pp [ \left| \{ i : \mu_i \leq T \} \right| = k ] &= \pp[ \left| \{ i : \lambda_i \leq T \} \right| + \left| \{ i : \lambda_i' \leq T \} \right| = 2k ]  \\
  &+ \pp[ \left| \{ i : \lambda_i \leq T \} \right| + \left| \{ i : \lambda_i' \leq T \} \right| = 2k +1] 
\end{align}
By the independence,
\begin{align}
\pp[ \left| \{ i : \lambda_i \leq T \} \right| +& \left| \{ i : \lambda_i' \leq T \} \right|  = n ] = \sum_{j=0}^n \pp[\left| \{ i : \lambda_i \leq T \} \right| =j] \times \pp[\left| \{ i : \lambda'_i \leq T \} \right| =(n-j) ]  \\
&= \sum_{j=0}^n \left( \pp[ \lambda_j \leq T ] - \pp [ \lambda_{j+1} \leq T ] \right) \times \left( \pp [ \lambda'_{n-j} \leq T ] - \pp [ \lambda'_{n-j+1} \leq T ] \right).
\end{align}
Taking the limit $N \to \infty$, we see
\begin{align}
\lim_{N \to \infty} \pp[ \left| \{ i : \lambda_i \leq T \} \right| &+ \left| \{ i : \lambda_i' \leq T \} \right|  = n ] \notag\\
=& \sum_{j=0}^n \left( \pp [ \chi_j \leq T ] - \pp [\chi_{j+1} \leq T ] \right) \times \left( \pp [ \chi'_{n-j} \leq T ] - \pp [ \chi'_{n-j+1} \leq T ] \right) \notag\\
=& \sum_{j=0}^n \pp [ \left| \{ i : \chi_i \leq T \} \right| = j ] \times \pp [ \left| \{ i : \chi'_i \leq T \} \right| = n-j ] \notag\\
=& \pp[ \left| \{ i : \chi_i \leq T \} \right|  + \left| \{ i : \chi'_i \leq T \} \right|  = n ]
\end{align}
This yields the claim. \qed

\bep
Let  $\{ \chi_i \}_i$ be  the Airy$_1$ random point field.  Then,
\beq
\left| \ee[ \left| \{ i : \chi_i \leq T \} \right| ] - \frac{2}{3} T^{3/2} \right| \leq C,
\eeq
and
\beq
\mathrm{Var} \left( \left| \{ i : \chi_i \leq T \} \right| \right) \leq C ( | \log (T)| + 1)
\eeq
\eep
\proof Let $\chi_i$ and $\chi'_i$ be two independent Airy$_1$ random point fields.  Let $Y$ be the random variable that is $k$ if there are $2k$ or $2k+1$ particles in the superposition $\{ \chi_i \}_i \cup \{ \chi'_i \}_i$ below $T$.  Then by the previous lemma, $Y$ has the same distribution as the number of particles in an Airy$_2$ random point field below $T$.  If $Z=\frac{1}{2} ( \left| \{ i: \chi_i \leq T \} \right| + \left| \{ i: \chi'_i \leq T \} \right|)-Y $, then $|Z| \leq C$.  The claim now follows. \qed

\subsection{Proof of Theorem \ref{thm:conv1}}
  Let us denote by $\N_T$ the random variable,
\beq
\N_T = \left| \{ i : \lambda_i \geq 2 - T N^{-2/3} \} \right|.
\eeq
Let $N^{4/15} \geq T \geq C_1$ where $C_1$ is the constant above.  We have by \eqref{eqn:vargoe},
\beq
\pp[ \lambda_k \geq 2 - T N^{-2/3} ] = \pp[ \N_T \geq k ] \leq  C \frac{ \log(T)}{ (\ee[ \N_T] - k)^2},
\eeq
as long as $k \geq \ee[ \N_T]$.  
Note that if $\gamma$ solves
\beq
\ee[ \N_\gamma] = x
\eeq
and $x \leq N^{ 2/5}$ then by \eqref{eqn:expgoe},
\beq
\gamma = \left( \frac{ 3 \pi x }{2} \right)^{2/3} + \O( (1+x^{1/3})^{-1} ).
\eeq
Assume that $k \leq N^{2/5}$.  Choosing now
\beq
 T = \left( \frac{3 \pi k}{2 } \right)^{2/3} - s
 \eeq 
 we see by \eqref{eqn:expgoe} that
 \beq
 k - \ee[ \N_T] \geq c s k^{1/3} - C,
 \eeq
 as long as
$s \leq (\frac{ 3 \pi k}{2})^{2/3} - C_1$.   Therefore, we have that
\beq 
\pp \left[ N^{2/3} (\lambda_k -2) \geq - \left( \frac{ 3 \pi k }{2} \right)^{2/3} + s \right] \leq C' \frac{ 1 + \log(k) }{ (s k^{1/3} - C)^2}.
\eeq
as long as $0 \leq s \leq  (\frac{ 3 \pi k}{2})^{2/3} - C_1$.  In particular, we see that there is a $K_1 >0$ so that for all $k \geq K_1$,
\beq \label{eqn:Gk}
\pp\left[ \bigcap_{N^{2/5  } \geq j \geq k } \left\{ N^{2/3} ( \lambda_j - 2 ) \leq -\left( \frac{ 3 \pi k }{2} \right)^{2/3} + \frac{1}{10} k^{2/3} \right\}  \right] \geq 1 - \frac{1}{k^{1/2}}.
\eeq
A similar argument gives
\beq
\pp \left[ N^{2/3} ( \lambda_k - 2 ) \leq - \left( \frac{ 3 \pi k }{2} \right)^{2/3}  - s \right] \leq C' \frac{ 1 + \log(k) + \log(1+s) }{( s k^{1/3} - C)^2}
\eeq
for $0 \leq s \leq N^{4/15}$ and $k \leq N^{2/5}$.  
From all of these estimates we find,
\beq
\ee \left[ \1_{ \{ N^{2/3} ( \lambda_k - 2 ) \leq -C_1 \} } \left| N^{2/3} ( \lambda_k - 2 ) + \left( \frac{ 3 \pi k }{2} \right)^{2/3}  \right| \right] \leq C \frac{ \log(k)^2}{k^{1/3}}, \qquad k \leq N^{2/5}.
\eeq
Denote by $\G_k$ the event on the left side of \eqref{eqn:Gk}.  Let $\eps >0$ and choose constants $C_\eps$ and $k_0$ so that
\beq
\pp[ \G_{k_0} ] \geq 1 - \eps, \qquad \pp[ | N^{2/3} ( \lambda_1 - 2 ) | \leq C_\eps ] \geq 1 - \eps.
\eeq
Let $\F$ denote the intersection of these two events.  Choose $K_2 \geq k_0$ so that $K_2 \geq 100 (C_\eps)^{3/2}$.   Fix also 
\beq
\delta_0 = \frac{1}{20}.
\eeq
By the choice of $K_2$ and the definition of $\F$ we have for any $k \geq K_2$ that on the event $\F$,
\beq
N^{2/3} ( \lambda_1 - \lambda_k ) \geq c k^{2/3}.
\eeq
 For any $k\geq K_2$ we then have,
\begin{align}
& \ee[ \1_F \left| \sum_{j = k}^{N^{\delta_0}}  \frac{1}{ N^{2/3} ( \lambda_1 - \lambda_j )}- \sum_{j=k}^{N^{\delta_0}} \int_{ (\frac{3 \pi (j-1)}{2} )^{2/3}}^{ ( \frac{3 \pi j }{2} )^{2/3} } \frac{1}{\pi \sqrt{x} } d x \right| ]  \notag\\
& \leq C \sum_{j=k}^{N^{\delta_0}} \frac{ C_\eps + C j^{-1/3} +  \ee[ \1_\F| N^{2/3} ( \lambda_j - 2 ) + \left( \frac{ 3 \pi j }{2} \right)^{2/3} |]}{j^{4/3}} \notag\\
& \leq C(C_\eps + 1 ) \log(k)^2 k^{-1/3}.
\end{align}
The outcome of all of this is that there for any $\eps >0$, there is a $K_3 >0$ so that for any $k \geq K_3$ we have the estimate,
\beq
\pp[ \left| \sum_{j = k}^{N^{\delta_0}}  \frac{1}{ N^{2/3} ( \lambda_1 - \lambda_j )}- \sum_{j=k}^{N^{\delta_0}} \int_{ (\frac{3 \pi (j-1)}{2} )^{2/3}}^{ ( \frac{3 \pi j }{2} )^{2/3} } \frac{1}{\pi \sqrt{x} } d x \right|  > \eps ] \leq \eps.
\eeq
By Theorem \ref{thm: rigidity},
\beq
\left| \frac{1}{N} \sum_{j = N^{\delta_0}}^N \frac{1}{ \lambda_1 - \lambda_j } - \int_{-2}^{\gamma_{N^{\delta_0}}} \frac{\rho_{\mathrm{sc}} (x) d x }{ 2 - x} \right| \leq N^{-c \delta_0-1/3}
\eeq
with overwhelming probability.  We write, for $k \geq K_3$,
\begin{align}
\frac{1}{N^{2/3} } \sum_{j=2}^N \frac{1}{ \lambda_1 - \lambda_j}  - N^{1/3}&= \frac{1}{N^{2/3}} \sum_{j=2}^{k} \frac{1}{ \lambda_1 - \lambda_j } - \int_{0}^{ ( \frac{ 3 \pi k}{2} )^{2/3} } \frac{ 1}{ \pi \sqrt{x}} \d x \\
&+ \frac{1}{N^{2/3}} \sum_{j=k+1}^{N^{\delta_0}} \frac{1}{ \lambda_1 - \lambda_j } - \int_{ ( \frac{ 3 \pi k }{2})^{2/3} }^{ ( \frac{ 3 \pi N^{ \delta_0 }}{2} )^{2/3} } \frac{1}{ \pi \sqrt{x}} \d x \\
&+ \frac{1}{N^{2/3}} \sum_{j=N^{\delta_0}+1}{N} \frac{1}{ \lambda_1 - \lambda_j } -  N^{1/3} \int_{-2}^{ \gamma_{N^{\delta_0}}} \frac{ \rho_{\mathrm{sc}} (x) }{ 2- x } d x \\
&+ \int_{0}^{ ( \frac{ 3 \pi N^{ \delta_0}}{2} )^{2/3}} \frac{1}{ \pi \sqrt{x}} d x - N^{1/3} \int_{\gamma_{N^{\delta_0}}}^{ 2} \frac{ \rho_{\mathrm{sc} (x)}}{ 2 - x } d x 
\end{align}
A calculation shows that the term on the last line is $\O ( N^{7\delta_0/6-2/3} ) = o (1)$ by our assumption on $\delta_0$.   Hence, for any bounded Lipschitz $F$ we see that for any $\eps >0$, there is a $k_1 = k_1 ( \eps)$ so that for any fixed $k > k_1$,
\begin{align}
&\limsup_{N \to \infty} \left| \ee \left[ F\left( N^{-2/3} \sum_{j=2}^N \frac{1}{ \lambda_1 - \lambda_j } - N^{1/3} \right) - F \left( N^{-2/3} \sum_{j=2}^k \frac{1}{ \lambda_1 - \lambda_j } -  \int_{0}^{ ( \frac{ 3 \pi k}{2} )^{2/3} } \frac{ 1}{ \pi \sqrt{x}} \right)\right] \right| \notag\\
 \leq &C ( ||F||_{Lip}) \eps.
\end{align}
With $\{ \chi_i \}_i$ denoting the particles of the Airy$_1$ process, much of the same calculations as above show that
\beq
\pp\left[ \chi_k \leq \left( \frac{ 3 \pi k }{ 2 } \right)^{2/3} - s \right] \leq \frac{C \log(k) }{ ( k^{1/3} s - C )^2}
\eeq
for $0 \leq s \leq k^{2/3} - C'_1$ some $C'_1 >0$, and
\beq
\pp \left[ \chi_k \geq \left( \frac{ 3 \pi k}{2} \right) ^{2/3} + s \right] \leq \frac{C \log(k) + C  \log (1+s)  }{ ( k^{1/3} s - C )^2}.
\eeq
Arguing as above, we see that for any $\eps >0$ there is an event $\F'$ with probability at least $1-\eps$ and a $K'_1 >0$ so that for all $k > K_1'$, we have
\beq
\ee[ \1_{\F'} \sum_{j \geq k } \left|  \frac{1}{ \chi_1 - \chi_j} +\int_{ ( \frac{3 \pi (j-1)}{2})^{2/3}}^{ ( \frac{ 3 \pi j }{2} )^{2/3} } \frac{ 1}{ \pi \sqrt{x} } d x \right| ] \leq C \log(k)^2 k^{-1/3}.
\eeq
From this we see that,
\beq
\limsup_{n \to \infty}  \sum_{j \geq n} \left|  \frac{1}{ \chi_1 - \chi_j} +\int_{ ( \frac{3 \pi (j-1)}{2})^{2/3}}^{ ( \frac{ 3 \pi j }{2} )^{2/3} } \frac{ 1}{ \pi \sqrt{x} } d x \right| = 0
\eeq
almost surely which proves that the limiting random variable $\Xi$ exists.  Moreover, we see that for any $\eps >0$ there is a $k_2$ so that for all $k > k_2$ and any bounded Lipschitz function $F:\mathbb{R}\rightarrow \mathbb{R}$ we have,
\begin{align}
&\left| \ee[ F ( \sum_{j=2}^\infty \frac{1}{ \chi_1 - \chi_j }  +\int_{ ( \frac{3 \pi (j-1)}{2})^{2/3}}^{ ( \frac{ 3 \pi j }{2} )^{2/3} } \frac{ 1}{ \pi \sqrt{x} } d x )] - \ee[ F ( \sum_{j=2}^k \frac{1}{ \chi_1 - \chi_j } +\int_{ ( \frac{3 \pi (j-1)}{2})^{2/3}}^{ ( \frac{ 3 \pi j }{2} )^{2/3} } \frac{ 1}{ \pi \sqrt{x} } d x ) ] \right| \notag\\
 \leq & C ||F||_{Lip} \eps.
\end{align}
On the other hand we know that for any finite $k$,
\beq
\frac{1}{N^{2/3}} \sum_{j=2}^k \frac{1}{ \lambda_1 - \lambda_j } \to  - \sum_{j=2}^k \frac{1}{ \chi_1 - \chi_j },
\eeq
where the convergence is in distribution as $N \to \infty$.  This yields the claim. \qed

\subsection{Proof of Theorem \ref{thm:mainconv}}

Let $H$ be the GOE matrix given by \eqref{eqn: H-def} and $M$ be the matrix \eqref{eqn:mdef}.  We choose a coupling so that $M_{ij} = H_{ij}$ for $i \neq j$.  
By Proposition \ref{prop:ev}, Theorem \ref{thm: rigidity} and Lemma \ref{lem:lr} there is an event of probability at least $1 - N^{-1/20}$ on which,
\beq
N^{1/3} \left| \frac{1}{N} \sum_{j=2}^N \frac{1}{ \lambda_1 (H) - \lambda_j (H) }- \frac{1}{ \lambda_1 (M) - \lambda_j (M) }\right| \leq N^{-c},
\eeq
for some $c>0$.  The result follows from this and Theorems \ref{thm1} and \ref{thm:conv1}. \qed

\section{Extension to $\langle | R_{12} |\rangle$} \label{sec:4m}
In this section we extend our results to the quantity $\langle|R_{12}| \rangle$.  Our starting point is the following elementary calculation,
\begin{align}
\langle |R_{12} | - q \rangle &= \left\langle \frac{ R_{12}^2-q^2}{|R_{12}| + q } \right\rangle \notag\\
&= \frac{1}{ 2 q } \langle R_{12}^2 - q^2 \rangle + \left\langle ( R_{12}^2 - q^2 ) \frac{ q - |R_{12} |}{(|R_{12} + q )2 q } \right\rangle
\end{align}
For the second term we have,
\begin{align}
\left| \left\langle ( R_{12}^2 - q^2 ) \frac{ |R_{12} |-q}{(|R_{12} + q )2 q } \right\rangle \right| = \left| \left\langle ( R_{12}^2 - q^2 )^2 \frac{1}{ (|R_{12}| + q )^2 2 q} \right\rangle \right| \leq \frac{1}{2 q^3} \langle (R_{12}^2 - q^2 )^2 \rangle.
\end{align}
Hence, if we can show that $\langle (R_{12}^2 - q^2 )^2 \rangle = o (N^{-1/3})$ with probability $1 - o (1)$, then the convergence of Theorem \ref{thm:mainconv} extends to $\langle |R_{12}| - q \rangle$.    

We expand,
\beq
\langle (R_{12}^2 - q^2)^2 \rangle = \langle R_{12}^4 \rangle - 2 q^2 \langle R_{12}^2 \rangle + q^4.
\eeq
We already calculated $\langle R_{12}^2 \rangle$ in Theorem \ref{thm:main} down to $o (N^{-2/3})$.  It remains to calculate the first term $\langle R_{12}^4 \rangle$.  The modification of the representation formula is,
\beq \label{eqn:4mrep}
\langle R_{12}^4 \rangle =  \frac{ \int_{\gamma- \i \infty}^{ \gamma + \i\infty}  \int_{\gamma- \i  \infty}^{ \gamma +\i \infty}  \e^{ \frac{N}{2} ( G(z) + G(w) ) }\left[6 \sum_{i=1}^N \frac{1}{ ( N \beta)^4 ( \lambda_i - w )^2 ( \lambda_i - z )^2 } +  3\left( \sum_{i=1}^N \frac{1}{ \beta^2 N^2 ( \lambda_i - z )(\lambda_i - w )} \right)^2  \right] \d  z \d w }{ \left( \int_{\gamma- \i \infty}^{\gamma+ \i \infty} \e^{ \frac{N}{2} G(z) } \d  z \right)^2}
\eeq
Since the function $G(z)$ appearing in the exponential is identical to what we encountered in considering $\langle R_{12}^2 \rangle$, the steepest descent analysis of the numerator is very similar to that in Section \ref{sec: sd}.  In light of this, we will use the same notation as in Section \ref{sec: sd}. 

Following along the argument of Section \ref{sec: sd} we see that, analogously to Lemma \ref{lem: hatcontours}, we can change the contour from the vertical line through $\gamma$ to $\hat{\Gamma}$ at an error exponential in $N^c$, for some $c>0$, on the event $\F_{\delta, \eps_1 }$ for sufficiently small $\eps_1$.

For the analog of Lemma \ref{lem: hattaylor} we similarly  derive the following estimates which all hold on the event $\F_{\delta, \eps_1}$ for $\eps_1 >0$ sufficiently small.  First,
\begin{align}
&\int_{\hat{\Gamma} \times \hat{\Gamma} } \e^{ \frac{N}{2} ( G ( z + \gamma )  + G ( w  + \gamma) -  2 G ( \gamma ) ) } \frac{1}{N ^4 ( z + \cb/N)^2( w + \cb / N )^2 }  \d z \d w \notag\\
= & \int_{\hat{\Gamma} \times \hat{\Gamma}} \e^{ \frac{N}{2} ( g ( z) + g ( w ) ) } \left( 1 + N z^2 \frac{ \tilm_N' ( \gamma ) }{4}  \right) \left( 1 + N w^2 \frac{ \tilm_N' ( \gamma ) }{4}  \right) \frac{1}{ N^4 ( z + \cb/N)^2 ( w + \cb / N )^2 } \notag\\
+ & \O \left( \frac{ N^{ 5 \kappa+ \eps_1 + 3 \delta}}{N^3} \right). \label{eqn:4m1}
\end{align}
Second,
\begin{align}
\int_{\hat{\Gamma} \times \hat{\Gamma}} \e^{ \frac{N}{2} ( G ( z + \gamma ) + G ( w + \gamma ) - 2 G ( \gamma ) )} \sum_{j=2}^N \frac{1}{N^4 (z+ \gamma - \lambda_j )^2 ( w+ \gamma - \lambda_j)^2 }  = \O \left( \frac{ N^{ 5 \kappa + \eps_1 + 3 \delta}}{N^3} \right).
\end{align}
Third,
\begin{align}
&\int_{\hat{\Gamma} \times \hat{ \Gamma}}  \e^{ \frac{N}{2} ( G ( z + \gamma ) + G ( w + \gamma ) - 2 G ( \gamma ) )}  \frac{1}{ N^2 ( z + \cb/N)(w + \cb/N) } \sum_{j=2}^N \frac{1}{N^2 ( z + \gamma - \lambda_j )( w + \gamma - \lambda_j ) }  \notag\\
= &\left(  \sum_{j=2}^N \frac{1}{N^2 ( \lambda_1 - \lambda_j )^2} \right) \notag\\
\times & \int_{\hat{\Gamma} \times \hat{\Gamma}} \e^{ \frac{N}{2} ( g (z) + g ( w ) ) }  \left( 1 + N z^2 \frac{ \tilm_N' ( \gamma ) }{4}  \right) \left( 1 + N w^2 \frac{ \tilm_N' ( \gamma ) }{4}  \right) \frac{1}{ N^2 ( z + \cb/N) ( w + \cb / N ) } \notag\\
+ & \O \left( \frac{ N^{ 5 \kappa+ \eps_1 + 3 \delta}}{N^3} \right) \label{eqn:4m3}
\end{align}
Finally,
\begin{align}
\int_{\hat{\Gamma} \times \hat{ \Gamma}}  \e^{ \frac{N}{2} ( G ( z + \gamma ) + G ( w + \gamma ) - 2 G ( \gamma ) )} \left( \sum_{j=2}^N \frac{1}{ N^2(z + \gamma - \lambda_j )( w + \gamma - \lambda_j )} \right)^2 = \O \left( \frac{ N^{ 5 \kappa+ \eps_1 + 3 \delta}}{N^3} \right).
\end{align}
In order to calculate the numerator of \eqref{eqn:4mrep} up to errors that are $o(N^{-2/3})$ we see that it suffices to compute the integrals in \eqref{eqn:4m1} and \eqref{eqn:4m3}.  We can proceed identically to Lemma  \ref{lem: contours-final} and pass to the rescaled variable $u$ being integrated over $\Gamma_r$ up to again an error exponential in $-N^{c}$ for some $c>0$.  The integral resulting from \eqref{eqn:4m3} is identical to \eqref{eqn:2mf2}, whereas the integral coming from \eqref{eqn:4m1} is 
\beq
\frac{1}{N} \int_{\Gamma_r}  \frac{ \e^{ ( \beta + \tilm_N ( \gamma ) ) u/2}}{ \sqrt{ 1 + u/\cb } } \left( 1 + \frac{ u^2 \tilm_N' ( \gamma ) }{ 4 N } \right) \frac{ d u }{ ( u + \cb )^2 }.
\eeq
In order to calculate this, we note the identities
\beq
\int_{\Gamma_{r, b} } \frac{ e^{ a z }}{ (z + b)^{5/2}} \d z = \e^{-ab} \frac{ 8  a^{3/2} \i}{3 } \sqrt{\pi}
\eeq
and
\beq
\int_{\Gamma_{r, b} } \frac{ e^{ a z } }{ (z + b)^{5/2}} z^2 \d z = \frac{ \e^{ - a b }}{ \sqrt{a}} \i \sqrt{\pi} \left( \frac{ 8 a^2 b^2}{3}  - 8 b a + 2 \right).
\eeq
From all of this, we see that we have derived the following estimate for the numerator of \eqref{eqn:4mrep}, with $a = (\beta + \tilm_N ( \gamma ) )/2$,
\begin{align}
& \int_{\gamma- \i \infty}^{ \gamma + \i\infty}  \int_{\gamma- \i  \infty}^{ \gamma +\i \infty}  \e^{ \frac{N}{2} ( G(z) + G(w) - 2 G ( \gamma) ) }  \notag\\
 \times & \left[6 \sum_{i=1}^N \frac{1}{ ( N \beta)^4 ( \lambda_i - w )^2 ( \lambda_i - z )^2 } +  3\left( \sum_{i=1}^N \frac{1}{ \beta^2 N^2 ( \lambda_i - z )(\lambda_i - w )} \right)^2  \right] \d  z \d w \notag\\
 = &  \left( \frac{ 2}{ \sqrt{a} } \frac{ \i \sqrt{ \pi } \e^{ - a \cb } \sqrt{ \cb }}{N} \right)^2  \left[ 16 a^4 -4 \frac{ \tilm_N' ( \gamma)  a^2 }{N}  + \frac{24 a^2}{N} \tilm_N' ( \lambda_1 ) \right] + \O \left( \frac{N^{5 \kappa + 3 \delta + \eps_1}}{N^3} \right)
\end{align}
whereas for the demoninator we have
\begin{align}
\left( \int_{\gamma- \i \infty}^{ \gamma+ \i \infty} \e^{ \frac{N}{2} G (z) } \right)^2 = \left( \frac{2}{ \sqrt{a}}  \frac{ \i \e^{ - a \cb} \sqrt{ \pi} \sqrt{ \cb }}{N} \right)^2 \left( 1 + \frac{ 3 \tilm_N' ( \lambda_1 ) }{4 N a^2 } \right) + \O \left( \frac{N^{5 \kappa + 3 \delta + \eps_1}}{N^3} \right).
\end{align}
From this and \eqref{eqn:r12final2} we see that
\begin{align}
\langle R_{12}^4 \rangle &= ( 1 - \beta^{-1} )^4 + 4 \frac{ ( \beta-1)^3}{ \beta^4} ( 1 + \tilm_N ( \lambda_1 ) ) + 6 \frac{ ( \beta-1)^2}{ \beta^4} ( 1 + \tilm_N ( \lambda_1 ) )^2 + 6 \frac{ ( \beta-1)^2}{\beta^4} \frac{ \tilm_N' ( \lambda_1 ) }{ N}  \notag\\
+& \O ( N^{-1+5 \kappa+\eps_1+3 \delta}),
\end{align}
and furthermore that
\begin{align}
\langle  ( R_{12}^2 - q^2 )^2 \rangle = \frac{ 8 ( \beta-1)^2}{ \beta^2} \frac{ \tilm_N'  ( \lambda_1 ) }{N} + 4 \frac{ ( \beta-1)^2}{\beta^4} ( 1 + \tilm_N ( \lambda_1 ) )^2 + \O (N^{-1+ 5\kappa+\eps_1+3 \delta}).
\end{align}
On the event $\F_{\delta, \eps_1}$ the first two terms are $\O( N^{-2/3+2 \delta+\eps_1} )$.

%

\appendix

\section{Zero-diagonal GOE} \label{a:diag}

Let $H$ be a GOE matrix as in \eqref{eqn: H-def}, and let $V$ be its diagonal and let
\beq
H = M + V
\eeq
so that $M$ is as in \eqref{eqn:mdef}. 
In this section we prove that with overwhelming probability, the extremal eigenvalues of $H$ and $M$ are close.
\bep \label{prop:ev} Let $\eps >0$.  
The following estimate holds for $i \leq N^{1/20}$ with overwhelming probability:
\beq
\left| \lambda_i (H) - \lambda_i (M) \right| \leq \frac{N^{\eps}}{N}.
\eeq
\eep
The proof follows from the Helffer-Sjoestrand formula, which we recall in \eqref{eqn: HS-formula}, and the following lemma providing control over the difference of Stieltjes transforms.
\bel \label{lem:stieltjes}
Denote by $m_M$ and $m_H$ the empirical Stieltjes transforms of $W$ and $H$.  Let $\eps >0$ and $\delta >0$.  With overwhelming probability, for any $N^{-\delta} \geq \eta \geq N^{\delta}/N$, and $|E| \leq 10$, we have
\beq
\left| m_M (z) - m_H (z) \right| \leq \frac{N^{\eps}}{N \eta} \left( \frac{1}{N \eta} + \Im [ \msc ] \right)
\eeq
\eel
The proof of the above lemma is based on the following resolvent expansion as well as two moment estimates which are the content of Lemmas \ref{lem:Ak} and \ref{lem:error} below.  We have,
\begin{align} \label{eqn:resolv}
\frac{1}{ M-z} =  \frac{1}{M-z} \sum_{k=1}^m  (V(M-z)^{-1} )^k+ \frac{1}{H-z} (V (M-z))^{-(m+1)}
\end{align}
Denote,
\beq
A_k := \frac{1}{M-z}  \left( V \frac{1}{M-z}  \right)^k, \qquad R(z) := \frac{1}{M-z}.
\eeq
Note that $G$ is independent of $V$.  We first prove,
\bel \label{lem:Ak}
Let $C>0$ be a constant.  On the event
\beq \label{eqn:Rass}
\max_{i, j} |R_{ij} | \leq C,
\eeq
we have for even $p$
\beq
\ee_V \left| \frac{1}{N} \tr A_k \right|^p \leq C(k,p) \left[ \frac{1}{ N \eta } \max_a \Im [ R_{aa} ] \right]^p,
\eeq
where $\ee_V$ denotes the expectation over $V$.
\eel
\remark This estimate is sub-optimal for $k \geq 2$ but we will not need a better estimate.  The next lemma below deals with the error term in the resolvent expansion. \qed

Before embarking on the proof, we record here the Ward identity, 
\beq \label{eqn:ward}
\sum_{a=1}^N  \left| \left( \frac{1}{A-z} \right)_{ab} \right|^2 = \frac{1}{ \eta} \Im \left( \frac{1}{A-z} \right)_{bb}
\eeq
for any self-adjoint matrix $A$.  This is a consequence of the spectral theorem (see Section 3 of \cite{KBG}).

\proof Denote by $\ulj = (j_1, \dots j_{kp} )$ and $\uli = ( i_1, \dots i_p)$ multi-indices in $kp$ and $p$ variables respectively, with $p$ even.   We  will use the $i$-indices to denote the summation coming from the trace, and the $j$-indices to be the summations coming from matrix multiplication.  Roughly, we are writing the trace as,
\beq
\frac{1}{N} \tr A_2 = \frac{1}{N} \sum_{i, j_1, j_2} R_{i j_1} V_{j_1} R_{ j_1 j_2 } V_{j_2} R_{j_2 i}.
\eeq
With this convention we then have,
\begin{align}
\ee_V \left| \frac{1}{N} \tr A_k \right|^p  = \frac{1}{N^p} \sum_{\uli} \sum_{\ulj} R^*_{i_1, j_1} R^*_{j_k, i_1} \dots R^*_{i_p j_{k(p-1)+1} } R^*_{i_p j_{kp} }  M( \ulj) \ee_V[ V_{j_1} \dots V_{j_{kp} } ].
\end{align}
where $R^*$ denotes either $R$ or $\bar{R}$ as appropriate, and $M( \ulj)$ is a monomial which contains all of the Green's function elements $R_{{j_a} j_{a+1}}$ that has indices only in  $\ulj$ (we separate out the matrix elements of $R$ that have an index in $\uli$).  The choice of $R$ or $\bar{R}$ or the form of $M$ will not be important for the calculation done here.  The $\ulj$ can be grouped into partitions of the $kp$ indices into coincidences.  That is,
\beq
\sum_{\ulj} = \sum_{ \P} \sum_{ \ulj \in \P}
\eeq
where the first sum is over partitions $\P$ on $kp$ elements, and the second summation means the sum over all $\ulj$ so that if $j_{a} = j_b$ whenever $a$ and $b$ are in the same block of $\P$ and $j_a \neq j_b$ whenever $a$ and $b$ are in distinct blocks of the partition.  The independence of the $V_j$ implies that unless the size of each block of the partion $\P$ is at least $2$, then the expectation vanishes.  Denote by $\P_2$ the set of such partitions.  Estimating $|M ( \ulj) | \leq C(k, p)$ (using the assumption \eqref{eqn:Rass}) we see that from this discussion,
\beq
\ee_V \left| \frac{1}{N} \tr A_k \right|^p  \leq \sum_{ \P \in \P_2 } \sum_{ \ulj \in \P } \frac{C}{N^{kp/2}} \frac{1}{N^p} \sum_{ \uli} |R_{i_1, j_1 } R_{i_1, j_k} \cdots R_{i_p, j_{kp}} |
\eeq
From the Ward identity \eqref{eqn:ward}, for any index $a$, $b$ we have
\beq
\frac{1}{N} \sum_{i_k=1}^N |R_{i_k a} R_{i_k b} | \leq \frac{1}{N} \sum_{i_k=1}^N |R_{i_k a } |^2+  |R_{i_k b }|^2 \leq \frac{1}{N \eta } \sup_k \Im [ R_{kk} ]
\eeq
and so
\beq
\frac{1}{N^p} \sum_{ \uli} |R_{i_1, j_1 } R_{i_1, j_k} \cdots R_{i_p, j_{kp}} | \leq \left( \frac{\sup_a \Im [ R_{aa} ]}{N \eta }  \right)^p.
\eeq
The summation over $\ulj \in \P$ for $\P \in \P_2$ has at most $N^{kp/2}$ terms, and $\P_2$ has cardinality bounded in terms of only $k$ and $p$, so we get the claim. \qed

\bel \label{lem:error}  Let $C>0$.  On the  event that $\max_{i , j} |R_{ij} | \leq C$ we have for even $p$ that,
\beq
\ee_V[ |(R (VR)^k )_{ab} |^p ] \leq C(k, p) \left( \frac{1}{N^{p(k/2-1)}} + \left( \max_{ i \neq j } |R_{ij} |^{p(k/2-1) } \right) \right)
\eeq
where $\ee_V$ denotes the expectation over $V$.
\eel
\proof We expand out the expectation similar to the proof of the above lemma.  We estimate $\sup_j |R_{aj}| \leq C$, and $\sup_j |R_{bj} | \leq C$ and obtain,
\beq
\ee_V[ |(R (VR)^k )_{ab} |^p ]  \leq \frac{C(k, p)}{N^{kp/2}} \sum_{ \P \in \P_2} \sum_{ \ulj \in \P } |M ( \ulj ) |
\eeq
where $\ulj$ is the following monomial in Green's function elements,
\beq
M ( \ulj ) = R_{j_1, j_2} R_{j_2, j_3 } \cdots R_{j_{k-1} j_k} R_{j_{k+1} j_{k+2}} \cdots R_{j_{kp-1} j_{kp}},
\eeq
i.e., it is the product of $R_{j_{i} j_{i+1}}$ except when $i=nk$ for any $n$.   Note that we have dropped any Green's function elements that involve the index $a$ or $b$, and  kept the ones involving only the $j_i$ indices.   We will use the estimate,
\beq
|M ( \ulj ) | \leq C(k, p ) \left( \max_{i \neq j } |R_{ij } | \right)^{ \mbox{\# off-diagonal}},
\eeq
and so we need to count how many off-diagonal entries appear in $M(\ulj)$ when $\ulj$ is in a specific partition $\P \in \P_2$.  Suppose that $\P \in \P_2$ has $\ell$ blocks.  Recall that $\ulj \in \P$ means that $j_a = j_b$ if and only if $a$ and $b$ are in the same block in the partition $\P$.
  Note that $M ( \ulj)$ contains $p(k-1)$ Green's function entries.  Denote the size of the $i$th block of $\P$ by $n_i \geq 2$.  There can be at most $n_i-1$ Green's function entries in the monomial $M ( \ulj)$ with whose indices $j_k$ and $j_{k+1}$ both appearin the $i$th block.  Therefore, there are at least
\beq
p(k-1) - \sum_{i=1}^\ell (n_i -1 )= pk - p - pk + \ell = \ell - p
\eeq
off-diagonal Green's function entries in the monomial $M ( \ulj)$.  Hence, for $\ulj \in \P$ where $\P$ has $\ell$ blocks,
\beq
| M ( \ulj ) | \leq  C(k, p) \left( \max_{i \neq j } |R_{ij} | \right)^{(\ell-p)_+}.
\eeq
The summation over $\ulj \in \P$ has less than $N^{\ell}$ terms, and so 
\begin{align}
\ee_V[ |(R (VR)^k )_{ab} |^p ]  &\leq C(k, p) \max_{ 1 \leq \ell \leq \frac{kp}{2} }  \frac{ N^\ell}{N^{kp/2}} \left( \max_{i \neq j } |R_{ij} | \right)^{(\ell-p)_+} \notag\\
&\leq C(k, p) \left( \frac{1}{N^{p(k/2-1)}} + \left( \max_{ i \neq j } |R_{ij} |^{p(k/2-1) } \right) \right)
\end{align}
This is the claim. \qed

\noindent{\bf Proof of Lemma \ref{lem:stieltjes}.} By the local semi-circle law, we have the estimates for $N^{\delta}/N \leq \eta \leq N^{-\delta}$,
\beq
\max_{i, j} |R_{ij} | \leq 2, \qquad \max_{i \neq j } |R_{ij} | \leq N^{-\delta/4}
\eeq
and
\beq
\max_a \Im [R_{aa} ] \leq N^{\eps} \left( \frac{1}{N \eta } + \Im [ m_{\mathrm{sc}} ] \right)
\eeq
with overwhelming probability.  Choose $m$ large enough so that $(m/2-1) \delta> 1000$.  Then by Lemma \ref{lem:error}, the final term in the resolvent expansion \eqref{eqn:resolv} is less than $N^{-100}$ with overwhelming probability, as it involves at most $N^2$ terms of the form $(R(VR)^m)_{ab}$ and $|(H-z)^{-1}_{ab}| \leq N$.  Finally, the other terms in the resolvent expansion are bounded using Lemma \ref{lem:Ak}. \qed

In order to prove Proposition \ref{prop:ev}, we will use the estimate on the Stieltjes transforms that we have just proved to find an estimate on traces of smoothed out indicator functions using the Helffer-Sj\"ostrand formula.  This is the content of the following lemma.
\bel \label{lem:hs}
Let $\eps_1, \eps >0$ and $\delta_f >0$ be arbitrary.   There is an event such that the following holds with overwhelming probability.  Suppose that $f$ is a smooth function so that $f=1$ on $[a, b]$ and $f=0$ outside of $[a - N^{\delta_f-1}, b + N^{\delta_f-1} ]$.  Assume that,
\beq
a \geq 2 - \kappa, \qquad |a-b| \leq \frac{ N^{\delta_f+\eps_1}}{N}.
\eeq
where $\kappa = N^{-1/2}$.  Assume that $\|f^{(k)}\|_{L^\infty} \leq C (N^{1-\delta_f})^k$ for $k=1, 2$.  Assume
\beq
0 < \delta_f < \frac{1}{20}, \qquad 0 < \eps_1 < \frac{1}{20}
\eeq
Then,
\beq
\left| \tr f (M) - \tr f (H) \right| \leq N^{\eps- \delta_f}.
\eeq
\eel
\proof  We work on the event that the estimate of Lemma \ref{lem:stieltjes} holds with $\eps >0$ and $\delta = \eps /10$.  We can assume that $100 \eps < \min \{ \eps_1, \delta_f \}$.  Fix,
\beq
\eta_1 = \frac{N^{\delta_1}}{N}, \qquad \eps < \delta_1 <\frac{1}{5}
\eeq
and let $\chi$ be a cut-off function so that $\chi (x) = 1$ for $|x| \leq \eta_1$ and $\chi (x) =0$ for $|x| > 2 \eta_1$ and $| \chi^{(k)} | \leq C(\eta_1)^{-k}$ for $k=1, 2$.  Denote
\beq
S = m_M - m_H
\eeq

Recall the Helffer-Sj\"ostrand formula (see, e.g., \cite{davies}): for $f\in C^2(\mathbb{R})$, define the almost analytic extension of $f$ by:
\[\tilde{f}(x+iy)=(f(x)+iyf'(x)) \chi (y) .\]
Then, for $x\in \mathbb{R}$, we have
\begin{equation}\label{eqn: HS-formula}
f(x)=\frac{1}{\pi}\int_{\mathbb{C}}\frac{\bar{\partial_z}\tilde{f}(z)}{x-z}\,\mathrm{d}z.
\end{equation}
Using this formula for we obtain,
\begin{align}
\left| \frac{1}{N} \tr f (M) - \frac{1}{N} \tr f (H) \right| &\leq  \left|  \int \int y f''(x) \chi (y) \Im [S(x + i y )] \d x \d y \right| \\
&+ \int \int |f(x) |  | \chi' (y) | | \Im S (x + \i y ) | \d x \d y \\
&+ \int \int |y f' (x) \chi' (y)  \Re S (x + \i y ) | \d x \d y.
\end{align}
Using Lemma \ref{lem:stieltjes} and the assumed $L^\infty$ bounds for $f'$ and $\chi'$, we see that the last term is bounded above by
\begin{align}
\int \int |y f' (x) \chi' (y)  \Re S (x + \i y ) | \d x \d y &\leq C \frac{N^\eps}{N} \left( \frac{1}{N \eta_1} + \sqrt{ \kappa + \eta_1 } \right) \notag\\
& \leq C N^{\eps-1}(N^{-\delta_1}+N^{-1/4} ),
\end{align}
where we used $\eta_1 < N^{-1/2}$ as well as $\Im [ \msc (z ) ] \leq C \sqrt{  | |E|-2| + \eta }$ (see, e.g., Section 3 of \cite{KBG}).  
Similarly, the second last term is bounded above by
\begin{align}
 \int \int |f(x) |  | \chi' (y) | | \Im S (x + \i y ) | \d x \d y &\leq C \|f\|_{L^1} \frac{1}{N \eta_1}  \left( \frac{1}{N \eta_1} + \sqrt{ \kappa + \eta_1 } \right) \notag\\
 &\leq  \frac{C N^{\delta_f+\eps_1}}{N} \left( N^{-2 \delta_1} + N^{-1/4} \right)
\end{align}
For the first term, we first estimate the contribution of $y \leq \eta_2  := N^{\eps}/N$.  From the fact that $y \to y \Im [m_M ( x + \i y ) ]$ is increasing (and the same for $m_H$) and the local semi-circle law, we find the estimate
\beq
\left| y \Im [ S ( x + \i y ) ]  \right| \leq 2\frac{N^\eps}{N} + C \eta_2 \Im [ m_{\mathrm{sc}} ]
\eeq 
which gives
\begin{align}
 \left|  \int \int_{|y| < \eta_2} y f''(x) \chi (y) \Im [S(x + i y )] \d x \d y \right| &\leq C\frac{N^{\eps} }{N^{\delta_f}} \left( \frac{N^{\eps}}{N} + C \eta_2 \sqrt{ \eta_2 + \kappa } \right) \notag\\
 &\leq C\frac{N^{2 \eps}}{N} N^{-\delta_f} 
\end{align}
For $\eta > \eta_2$ we get,
\begin{align}
\left|  \int \int_{|y| > \eta_2} y f''(x) \chi (y) \Im [S(x + i y )] \d x \d y \right| &\leq C N^{\eps} \frac{N}{N^{\delta_f}}  \int_{\eta_2 < y <2 \eta_1} \left( \frac{1}{N^2 y} + \frac{1}{N} \sqrt{y+\kappa}  \right) \notag \\
&\leq C \left( \frac{N^{2 \eps}}{N^{1+\delta_f}} + N^{\eps-\delta_f} \sqrt{\kappa} \eta_1 + N^{\eps-\delta_f} \eta_1^{3/2} \right) \leq C N^{2 \eps-1-\delta_f}
\end{align}
 We have proven,
\beq
\left|  \tr f (M) -  \tr f (H) \right| \leq N^{3 \eps} \left( N^{-\delta_1} + N^{\delta_f+\eps_1-2 \delta_1 } +N^{-\delta_f } \right).
\eeq
The claim follows from choosing $\delta_1 = \delta_f + \eps_1 + \frac{1}{10}$. \qed

\noindent{\bf Proof of Proposition \ref{prop:ev}.}
Let now $0 < \delta_f < \frac{1}{20}$ and $\eps >0$.  Let us denote by $\lambda_i$ the eigenvalues of $H$ and by $\mu_i$ the eigenvalues of $M$.  Applying Lemma \ref{lem:hs} to $f$ with $|a-b| = N^{\delta_f}/N$, and $a> \lambda_1 + N^{\delta_f-1}$,  we see first that
\beq
\mu_1 \leq \lambda_1 + N^{\delta_f-1}.
\eeq
(or else for some $a > \lambda_1 + N^{\delta_f-1}$ we would have $\tr f (M) \geq 1$, contradicting the estimate proven in Lemma \ref{lem:hs} and the fact that $\tr f (H) = 0$ for such $f$).  
Reversing the roles of $\mu_1$ and $\lambda_1$ we then get that 
\beq
| \mu_1 - \lambda_1 | \leq N^{\delta_f-1}.
\eeq

Let $k_1$ be the smallest index so that
\beq
\lambda_{k_1} - \lambda_{k_1+1} > 10 \frac{N^{\delta_f}}{N}
\eeq
and let $J_1 = [[1, k_1]]$.  Similarly, let $k_2 > k_1$ be the smallest index so that
\beq
\lambda_{k_2} - \lambda_{k_2+1} > 10 \frac{N^{\delta_f}}{N},
\eeq
and let $J_2 = [[k_1+1, k_2]]$, and define $J_i = [[k_{i-1}+1, k_i]]$ and so on.  Let $\ell$ be the smallest integer so that
\beq
[[1, N^{1/20}]] \subseteq J_1 \cup J_2 \cup \cdots \cup J_\ell.
\eeq
By rigidity we have that for $j - i > N^{\eps}$,
\beq
|\lambda_i - \lambda_j | > c \frac{ j-i}{ N^{2/3} ( j^{1/3} )} \gg N^{\delta_f} \frac{j-i}{N},
\eeq
if $j \leq N^{1/20}$.  Therefore,
\beq
|J_i | \leq N^{\eps}, \mbox{ for } i \leq \ell.
\eeq
First, using Lemma \ref{lem:hs} we see that there are no eigenvalues $\mu_i$ in the interval
\beq
[\lambda_{k_i+1}+ N^{\delta_f-1}, \lambda_{k_i}-  N^{\delta_f-1} ],
\eeq
for $i \leq \ell$, by taking $ a > \lambda_{k_i+1} + N^{\delta_f-1}$ and $b < \lambda_{k_i} - N^{\delta_f-1}$, and $|a-b| = N^{\delta_f-1}$.  Next, we apply Lemma \ref{lem:hs} with the choice $b= \lambda_{k_i+1} +  N^{\delta_f-1}$ and $a = \lambda_{k_{i+1} } -  N^{\delta_f-1}$.  Note that since $|J_i | \leq N^{\eps}$, we see that the length of $[a, b]$ in this case is less than $C N^{\eps+\delta_f-1}$, and so the lemma applies, which gives
\beq
\tr f ( M) = |J_i | + o (1).
\eeq
Since we have already shown that there are no eigenvalues $\mu_i$ in the intervals $[\lambda_{k_{i+1}}-2N^{\delta_f-1} , \lambda_{k_{i+1}}- N^{\delta_f-1} ]$ and $[\lambda_{k_i+1}+N^{\delta_f-1}, \lambda_{k_i+1} + 2 N^{\delta_f-1} ]$, it follows that the quantity $\tr f (M)$ is precisely the number of eigenvalues in the interval $[\lambda_{k_{i+1} } - N^{\delta_f-1} , \lambda_{k_i+1} + N^{\delta_f-1} ]$.  This must be an integer, and so it equals $|J_i|$.   The claim follows.  \qed

\vspace{5 pt}
\noindent\textbf{Acknowledgements}. We wish to thank Jinho Baik for sharing some computations from work in preparation, and Hao Wu for pointing out an error in an earlier draft of the manuscript.  P.S. thanks Vu Lan Nguyen for discussions on the spherical SK model at low temperature.  B.L. thanks Amol Aggarwal for many illuminating and useful discussions.

\bibliography{ssk_bib}
\bibliographystyle{abbrv}

\end{document}